\documentclass[11pt]{articulofederico}

\usepackage{amssymb}
\usepackage{dcolumn}
\usepackage{enumerate}
\usepackage{subfigure}
\usepackage{epsfig}

\usepackage{amsfonts} 

\def\NantelYoung#1{\vbox{\smallskip\offinterlineskip
    \halign{&\vbox{##}\kern-\Thickness\cr #1}}}

\newdimen\Squaresize \Squaresize=15pt
\newdimen\Thickness \Thickness=.5pt
\newdimen\Correction \Correction=7pt

\def\Vide#1{\hbox{
       \vbox to \Squaresize{\vss
          \hbox to \Squaresize{\hss#1 \hss}\vss}
    \hskip-\Correction}
   \kern-\Thickness}

\def\Carre#1{\hbox{\vrule width \Thickness
   \vbox to \Squaresize{\hrule height \Thickness\vss
      \hbox to \Squaresize{\hss$\scriptstyle#1$\hss}
   \vss\hrule height\Thickness}
   \unskip\vrule width \Thickness}
   \kern-\Thickness}

\usepackage{amssymb}
\usepackage{enumerate}
\usepackage{graphicx}

\usepackage[latin1]{inputenc}
\usepackage{amsmath,amssymb,enumerate,amsthm,flafter}

\newtheorem{thm}{Teorema}[section]
\newtheorem{prop}[thm]{Proposici\'on}
\newtheorem{cor}[thm]{Corolario}

\newtheorem{obs}[thm]{Observaci\'on} 
 
\theoremstyle{definition}
 \newtheorem{defn}{Definici\'on}[section]
 
 \newtheorem{question}[defn]{Pregunta}

\newenvironment{dem}{{\noindent {\bf \it Demostraci\'on: }}}%
{\hspace*{\fill}\nolinebreak[1]\hspace*{\fill}$\Box$\vspace{.1 in}}%
\renewcommand{\figurename}{Figura}

\newcommand{\notI}{\overline{I}}

\newcommand{\notJ}{\overline{J}}

\newcommand{\ssm}{\smallsetminus}

\renewcommand{\refname}{\textsf{Referencias}}
\renewcommand{\abstractname}{\textsf{Resumen}}

\title{\textsf{Tres lecciones en combinatoria algebraica.  \\
\normalsize{I. Matrices
totalmente no negativas y funciones sim\'etricas.}}}

 \author{\textsf{
 Federico Ardila\footnote{\textsf{San Francisco State University, San Francisco, CA, USA y Universidad de Los Andes, Bogot\'a, Colombia, federico@sfsu.edu -- financiado por la CAREER Award DMS-0956178 y la beca DMS-0801075 de la National Science Foundation de los Estados Unidos, y por la SFSU-Colombia Combinatorics Initiative.}}}\qquad
\textsf{Emerson Le\'on\footnote{\textsf{Freie Universit\"at Berlin, Alemania,  emerson@zedat.fu-berlin.de -financiado por el Berlin Mathematical School.}}}\\
\textsf{Mercedes Rosas\footnote{\textsf{Universidad de Sevilla, Espa\~na, mrosas@us.es -- financiada por los proyectos MTM2007--64509 del Ministerio de Ciencias e Innovaci\'on de Espa\~na y FQM333 de la Junta de Andalucia.}}}\qquad
\textsf{Mark Skandera\footnote{\textsf{Lehigh University, 
Bethlehem, PA, USA, 
mas906@math.lehigh.edu -- financiado por la beca H98230-11-1-0192 de la National Security Agency de los Estados Unidos.}}}
}

\date{}

\begin{document}

\maketitle

\begin{abstract} 
En esta serie de tres art\'\i culos, damos una exposici\'on de varios resultados y problemas abiertos en tres \'areas de la combinatoria algebraica y geom\'etrica:  las matrices totalmente no negativas,
las representaciones
de $S_n$, y los arreglos de hiperplanos. Esta primera parte presenta una introducci\;on a las matrices totalmente no negativas, y su relaci\'on con las funciones sim\'etricas.
\end{abstract}



\section{\textsf{Introducci\'on}\label{s:intro}}

En marzo de 2003 se llev\'o a cabo el Primer Encuentro Colombiano de Combinatoria en Bogot\'a, Colombia. Como parte del encuentro, se organizaron tres minicursos, dictados por Federico Ardila, Mercedes Rosas, y Mark Skandera. Esta serie resume el material presentado en estos cursos en tres art\'{\i}culos:   \emph{I. Matrices
totalmente no negativas y funciones sim\'etricas} \cite{ALRS1}, 
\emph{II. Las funciones sim\'etricas y la teor\'{\i}a de las representaciones} \cite{ALRS2}, y \emph{III. Arreglos de hiperplanos.} \cite{ALRS3} 

Uno de los mensajes que quisimos transmitir en nuestros minucursos es la idea de que muchos problemas en matem\'aticas se entienden mejor mediante el estudio de los {\em objetos combinatorios} que los subyacen, como grafos, particiones y conjuntos parcialmente ordenados (posets). En este primer art\'\i culo, veremos la utilidad de estos objetos combinatorios en el estudio de las matrices totalmente no negativas y las funciones sim\'etricas. 

\section{\textsf{Matrices totalmente no negativas}\label{s:mtnn}}
Sea $A$ una matriz
$n \times n$ 
y sean $I$ e $I'$ dos subconjuntos de
$[n] = \{1,\dotsc,n\}$.
Definimos 
$A_{I,I'}$ como la submatriz de $A$ que se construye a partir de $A$
con las entradas correspondientes a las filas en
$I$ y las columnas en $I'$.  El determinante de esta submatriz
recibe el nombre del {\em menor $(I,I')$ de $A$}, 
y se escribe 
\begin{equation*}
\Delta_{I,I'} = \det A_{I,I'}.
\end{equation*}
(Para hablar del menor $(I,I')$ necesitamos que $|I|=|I'|$.)
Decimos que una matriz es {\em totalmente no negativa} si todos sus
menores son mayores que o iguales a cero y {\em totalmente positiva}
si todos los menores son mayores que cero.

Un ejemplo de una matriz totalmente no negativa es
\begin{equation}\label{eq:matrizdecaminos}
\begin{bmatrix}
5 & 6 & 3 & 0 \\ 
4 & 7 & 4 & 0 \\
1 & 4 & 4 & 2 \\
0 & 1 & 2 & 3 
\end{bmatrix}.
\end{equation} 
Se puede verificar, aunque sea tedioso, que todos los menores de la matriz 
anterior, como por ejemplo
\begin{equation*}
\Delta_{\{1,2\},\{1,3\}} = \det 
\begin{bmatrix}
 5 & 3 \\
 4 & 4
\end{bmatrix}
= 8,
\end{equation*}
son mayores que o iguales a cero.

Las matrices totalmente no negativas aparecen en varias \'areas de la
matem\'atica como ecuaciones diferenciales~\cite{GantKreinOsc},
ra\'{\i}ces de polinomios~\cite{ASW}, \cite[p.\,241]{Ando}, 
\nocite{Gant1} 
\cite[Cap.\,XV]{Gant2},
procesos estoc\'asticos~\cite{KMG}, \cite{FominLoop},
matroides~\cite{LinVrep},
grupos de Lie~\cite{LusztigTP}, 
geometr\'\i a algebraica~\cite{FominDouble}, \cite{FominShapStrat},
\cite{PostQuantumBruhat},
y conjuntos parcialmente ordenados~\cite{Skan3+1}, \cite{SkanReed}.
Estas matrices aparecen tambi\'en en 
econom\'\i a~\cite{KarlinGame}, 
ingenier\'\i a el\'ectrica~\cite{CurtisCircular} y
qu\'\i mica~\cite{SachsHex}.

En los art\'\i culos de 
Karlin y McGregor~\cite{KMG} y Lindstr\"om~\cite{LinVrep},
estas matrices surgieron en conexi\'on con ciertos grafos
denominados {\em redes planas}.
Una {\em red plana de orden $n$} es un grafo dirigido ac\'\i clico
que se puede dibujar en el plano sin que las aristas se corten en
puntos diferente de los v\'ertices, tal que hay $2n$ v\'ertices
distinguidos: $n$ {\em fuentes} (a las cuales no llega
ninguna arista) y $n$ {\em destinos} (de los cuales no sale ninguna
arista).

\begin{figure}[h]
\centerline{\psfig{figure=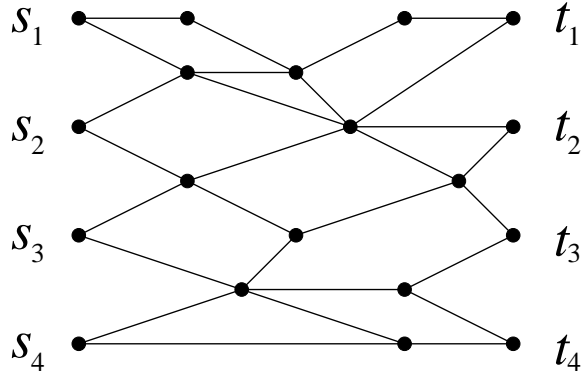,width=.50\textwidth}}
\caption{Una red plana de orden $4$.}
\label{f:nplanar}
\end{figure}

La Figura~\ref{f:nplanar} muestra una red plana de orden $4$.
En nuestras figuras, siempre ubicaremos a las fuentes y a los destinos 
en los extremos
izquierdo y derecho del grafo, respectivamente, y les daremos las etiquetas 
$s_1,\dotsc,s_n, t_n,\dotsc,t_1$ de arriba hacia abajo. 
Tambi\'en daremos por supuesto que cada arista
est\'a orientada hacia la derecha (a menos que indiquemos
otra orientaci\'on).

Dada una red plana, definimos la {\em matriz de caminos} 
$A = [a_{ij}]$ de $G$ como
\begin{equation*}
a_{ij} = \# \text{ caminos desde $s_i$ hasta $t_j$ en $G$.}
\end{equation*}
La matriz de caminos de la red plana en la Figura~\ref{f:nplanar} se muestra en (\ref{eq:matrizdecaminos}).

\begin{thm}\label{t:linlem}
(Lema de Lindstr\"om)
La matriz de caminos de una red plana
es totalmente no negativa.  
El menor $\Delta_{I,J}$ es igual al n\'umero de familias 
de $n$ caminos disjuntos desde las fuentes $I$ hasta los destinos $J$.
\end{thm}

\begin{dem}
Usando inducci\'on, supongamos que para todas las redes planas
de orden menor que $n$, las matrices de caminos correspondientes
son totalmente no negativas.  (Esto es cierto para $n=1$.)
Sea $G$ una red plana de orden $n$ con matriz de caminos $A$.
Dado que cada submatriz de $A$ es la matriz de caminos de una subred
de $G$, es suficiente mostrar simplemente que el determinante de $A$
es igual al n\'umero de familias
$\pi = (\pi_1,\dotsc, \pi_n)$ de $n$ caminos en $G$ en las cuales
$\pi_i$ es un camino desde $s_i$ hasta $t_i$, por $i=1,\dotsc,n$,
y en las cuales no hay ning\'un v\'ertice que pertenece a dos caminos.

Veamos el producto $a_{1,\sigma(1)} \cdots a_{n,\sigma(n)}$, que aparece en el determinante
de $A$,
\begin{equation*}
\det A = \sum_{\sigma \in S_n}
\text{sgn}(\sigma)a_{1,\sigma(1)} \cdots a_{n,\sigma(n)}.
\end{equation*}
Este producto es igual al n\'umero de familias de $n$ caminos en $G$
en las cuales el camino $i$ va desde $s_i$ hasta $t_{\sigma(i)}$, donde se permiten intersecciones.
Sea $\pi = (\pi_1,\dotsc,\pi_n)$
una familia en la cual hay por lo menos una intersecci\'on; sea $j$ el menor \'{\i}ndice tal que los caminos $\pi_j$ y $\pi_{j+1}$ se intersectan, y consideremos su primera intersecci\'on.
El conjunto de aristas utilizadas por $\pi$ es 
el mismo conjunto de aristas utilizadas por la familia
\begin{equation*}
\pi^* = (\pi_1,\dotsc,\pi_{j-1}, \pi_j', \pi_{j+1}',\pi_{j+2},\dotsc,\pi_n)
\end{equation*} 
que se construye intercambiando a los caminos $\pi_j$ y $\pi_{j+1}$ a partir
de su primera intersecci\'on.
(Vea la Figura~\ref{f:swap}.)

\begin{figure}[h]
\centerline{\psfig{figure=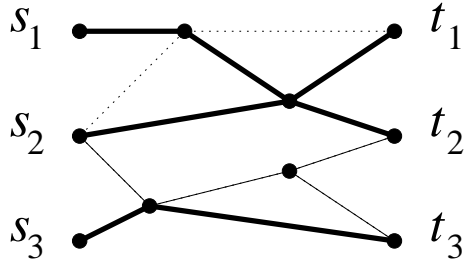,width=.40\textwidth}}
\caption{Un conjunto de aristas que define una familia de caminos desde
$s_1,s_2,s_3$ hasta $t_1,t_2,t_3$ (respectivamente) 
y una familia de caminos desde
$s_1,s_2,s_3$ hasta $t_2,t_1,t_3$ (respectivamente).}
\label{f:swap}
\end{figure}

La familia $\pi^*$ es contada por 
$a_{1,\sigma(1)} \cdots a_{j-1,\sigma(j-1)} a_{j,\sigma(j+1)} a_{j+1,\sigma(j)} a_{j+2,\sigma(j+2)} \cdots a_{n,\sigma(n)}$,
y aparece en $\det A$ con el signo contrario al de $\pi$, ya que las dos permutaciones difieren en una transposici\'on. Tenemos entonces que la contribuci\'on total de estas dos familias al 
determinante es cero. Observemos adem\'as que $(\pi^*)^*=\pi$. 

Por lo tanto, en la expresi\'on del determinante, las familias $\pi$ que se intersectan est\'an emparejadas de manera que la contribuci\'on total de cada pareja es $0$. Entre tanto, las familias $\pi$ que no se intersectan tienen que unir a $s_i$ con $t_i$ para cada $i$; por lo tanto, en la expresi\'on del determinante, todas aparecen con signo positivo. El resultado se sigue. \end{dem}

\begin{figure}[h]
  \begin{center}
    \mbox{
      \subfigure{\psfig{figure=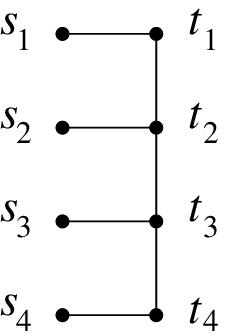,width=.20\textwidth}
      \label{f:single}} \qquad
      \subfigure{\psfig{figure=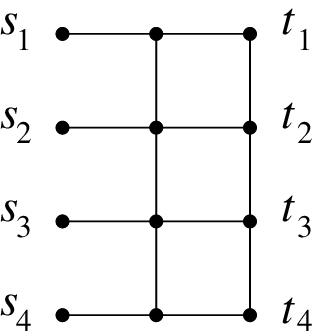,width=.285\textwidth}
      \label{f:double}}
      }
    \caption{Dos redes planas.}
    \label{f:singledouble}
  \end{center}
\end{figure}

El Lema de Lindstr\"om nos da una demostraci\'on f\'acil de que ciertas
matrices son totalmente no negativas.
Por ejemplo las matrices
\begin{equation}\label{eq:singledouble}
A = \begin{bmatrix}
1 & 0 & 0 & 0 \\ 
1 & 1 & 0 & 0 \\
1 & 1 & 1 & 0 \\
1 & 1 & 1 & 1
\end{bmatrix}, 
\qquad
B = \begin{bmatrix}
1 & 0 & 0 & 0 \\ 
2 & 1 & 0 & 0 \\
3 & 2 & 1 & 0 \\
4 & 3 & 2 & 1
\end{bmatrix}
\end{equation}
son totalmente no negativas porque son las matrices de caminos de las redes
planas de la Figura~\ref{f:singledouble}, donde las aristas tienen orientaci\'on
hacia la derecha y hacia arriba.

Lindstr\"om~\cite{LinVrep} (y Karlin y McGregor~\cite{KMG}) demostraron el
Teorema~\ref{t:linlem} de forma más general.
Dado una red plana $G$ cuyas aristas tienen pesos positivos y reales,
definimos el {\em peso} de un camino como el producto de los pesos de
sus aristas y definimos la {\em matriz de pesos} 
$A = [a_{ij}]$ de $G$ como
\begin{equation*}
a_{ij} = \text{suma de pesos de caminos desde $s_i$ hasta $t_j$ en $G$.}
\end{equation*}
La matriz de pesos tambi\'en es totalmente no negativa y sus menores
se pueden interpretar de manera similar:

\begin{thm}\label{t:wlinlem}
La matriz de pesos de una red plana es totalmente no negativa.
El menor $\Delta_{I,J}$ es igual a la suma 
de los pesos de las familias de caminos disjuntos desde las fuentes $I$ hasta 
los destinos $J$.
\end{thm}

Observe que el Teorema~\ref{t:linlem} es el caso especial del 
Teorema~\ref{t:wlinlem} que corresponde a darle peso $1$ a cada arista.
La demostraci\'on del Teorema~\ref{t:wlinlem} es muy similar a la del
Teorema~\ref{t:linlem}.

La Figura~\ref{f:singledouble} ilustra una propiedad interesante de las matrices totalmente no negativas. Las matrices en (\ref{eq:singledouble}) satisfacen la ecuaci\'on
$A^2 = B$, y se puede construir una red plana que corres\-ponde a $B$
utilizando dos copias de la red que corresponde a $A$.
En general, concatenar redes planas corresponde a multiplicar matrices totalmente no negativas. Este hecho se sigue directamente de la definici\'on de multiplicaci\'on
de matrices. 

\begin{figure}[h]
\centering
\mbox{
      \subfigure{\epsfig{figure=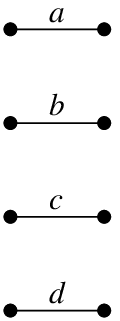,width=.1\textwidth}
      \label{f:tnnfactor1}} \qquad\qquad
      \subfigure{\epsfig{figure=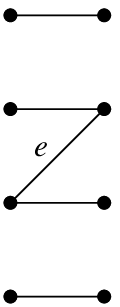,width=.1\textwidth}
      \label{f:tnnfactor2}} \qquad\qquad
      \subfigure{\epsfig{figure=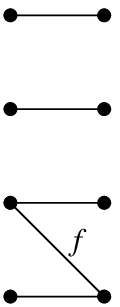,width=.1\textwidth}
      \label{f:tnnfactor3}}
      }
    \caption{Tres redes planas elementales.}
    \label{f:tnnfactors}
\end{figure}

La correspondencia entre multiplicaci\'on de matrices y concatenaci\'on de redes tiene una
consecuencia extraordinaria.
A. Whitney~\cite{WReduction} y Loewner~\cite{Loewner} demostraron que
cada matriz totalmente no negativa e invertible se puede factorizar como
un producto de matrices totalmente no negativas
\begin{equation}\label{eq:factortnn}
L_1 \cdots L_m D U_1 \cdots U_m,
\end{equation}
tal que  $D$ es una matriz diagonal,
$L_i$ tiene la forma $I + cE_{j+1,j}$ y
$U_i$ tiene la forma $I + cE_{j,j+1}$,
donde $E_{k,\ell}$ es la matriz que tiene un $1$ en la posici\'on $k,\ell$ 
y ceros en todas las dem\'as posiciones.
Cryer~\cite{CryerProp} extendi\'o este resultado a {\em cualquier} matriz totalmente no negativa.
Est\'a claro que los factores que aparecen en (\ref{eq:factortnn})
son matrices de
pesos de redes planas~\cite{BrentiCTP}; 
la Figura~\ref{f:tnnfactors} muestra las redes planas cuyas matrices de pesos
son
\begin{equation*}
\begin{bmatrix}
a & 0 & 0 & 0 \\ 
0 & b & 0 & 0 \\
0 & 0 & c & 0 \\
0 & 0 & 0 & d
\end{bmatrix}, 
\qquad
\begin{bmatrix}
1 & 0 & 0 & 0 \\ 
0 & 1 & 0 & 0 \\
0 & e & 1 & 0 \\
0 & 0 & 0 & 1
\end{bmatrix}, 
\qquad
\begin{bmatrix}
1 & 0 & 0 & 0 \\ 
0 & 1 & 0 & 0 \\
0 & 0 & 1 & f \\
0 & 0 & 0 & 1
\end{bmatrix}. 
\end{equation*}
Vemos entonces que las matrices de pesos no son solamente una clase importante de matrices totalmente no negativas; las contiene a todas. 

\begin{thm}\label{t:tnnconverse}
Cada matriz $n \times n$ que es totalmente no negativa es 
la matriz de pesos de una red plana de orden $n$.
\end{thm}

\begin{figure}[h]
\centerline{\psfig{figure=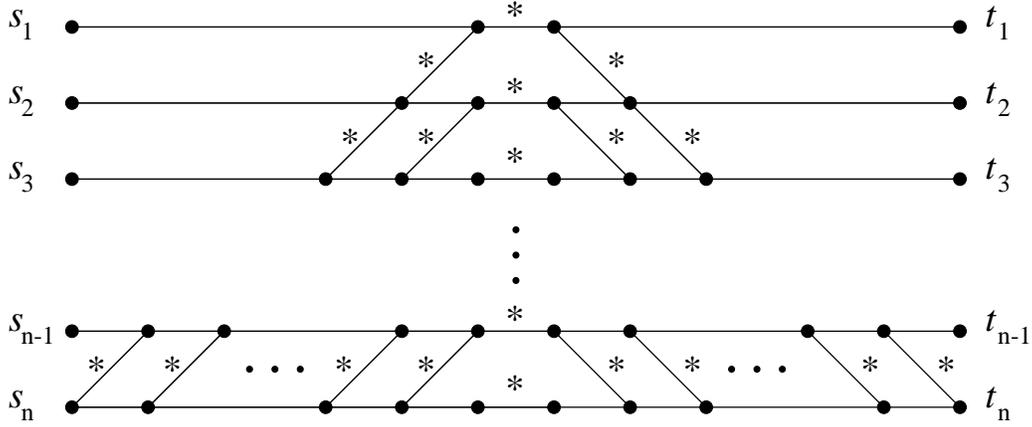,width=.90\textwidth}}
\caption{Una red plana 
de orden $n$ 
que corresponde a una matriz $n \times n$ totalmente positiva e invertible.}
\label{f:tnnfoz}
\end{figure}

Aunque no existe una {\em \'unica} red plana que corresponda a
una fija matriz totalmente no negativa, 
si la matriz es invertible, podemos elegir una red plana
que tenga la forma mostrada en la Figura~\ref{f:tnnfoz}~\cite{FominTPTest}.
Entonces podemos calcular de manera \'unica los pesos indicados por $*$.  
Estos pesos son funciones racionales de las entradas 
de la matriz.  Entonces si la matriz tiene
entradas enteras, los pesos ser\'an n\'umeros racionales.

\begin{obs}\label{o:fomin}
Para cada matriz $A$ totalmente no negativa con entradas
enteras, existe un entero $k$ tal que $kA$ es la matriz de caminos
de una red plana sin pesos.
\end{obs}

Como es de esperarse, el hecho de que una matriz sea totalmente no negativa
nos da informaci\'on sobre sus
valores propios.
\begin{thm}\label{t:tnneigen}
Todos los valores propios de una matriz totalmente no negativa son n\'umeros
reales.
\end{thm}

\begin{dem}
Una demostraci\'on sencilla utiliza \'algebra exterior.
(Vea \cite[pp.\,167-172]{Ando}.)
Dada una matriz $A$ de $n \times n$,
la potencia exterior $k$ de $A$, escrita $\wedge^k A$,
se puede definir como la matriz
$\binom{n}{k} \times \binom{n}{k}$ cuyas filas y columnas son los subconjuntos
de $[n]$ con $k$ elementos y cuya entrada $I,J$ es
$\Delta_{I,J}$, el menor $(I,J)$ de $A$.
Los valores propios de $\wedge^k A$ son todos los $\binom{n}{k}$ 
productos de  $k$ valores propios de $A$  (contando
multiplicidades).

Veamos primero el caso de una matriz totalmente positiva.
Sea $A$ una matriz totalmente positiva cuyos valores propios son
$\lambda_1,\dotsc,\lambda_n$, con
\begin{equation*}
|\lambda_1| \geq \cdots \geq |\lambda_n|.
\end{equation*}
Las entradas de $A$, al ser los menores de $1 \times 1$, son positivas.
Entonces el Teorema de Perron-Frobenius
(Vea \cite[p.189]{Luen}.)~nos dice que $\lambda_1$ es el \'unico valor propio
de m\'aximo valor absoluto y que $\lambda_1$ es positivo.
Aplicando el Teorema de Perron-Frobenius a 
$\wedge^2 A, \wedge^3 A,\dotsc, \wedge^n A$ vemos que cada
matriz tiene un \'unico valor propio
de valor absoluto m\'aximo.  Estos valores propios son
$\lambda_1\lambda_2, \lambda_1\lambda_2\lambda_3, \dotsc,
\lambda_1\cdots\lambda_n$, respectivamente.
Dado que cada valor propio es positivo, vemos entonces
que los n\'umeros
$\lambda_1,\dotsc,\lambda_n$ tambi\'en son positivos.

Ahora supongamos que $A$ es totalmente no negativa.  
Un resultado  conocido dice que $A$ es el l\'\i mite de una sucesi\'on
de matrices totalmente positivas.
(Vea~\cite[Tm.\,2.7]{Ando} o ~\cite[Tm.\,11]{FominTPTest}.)
Dado que los valores propios de estas matrices 
convergen a los valores propios de $A$, \'estos tienen que ser reales
y no negativos.
\end{dem}

\begin{figure}[h]
\centerline{\psfig{figure=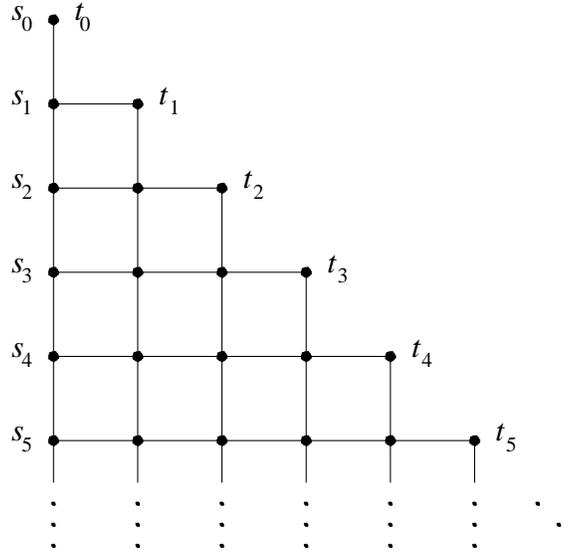,width=.50\textwidth}}
\caption{La red plana que corresponde a la matriz infinita de coeficientes
binomialeas.}
\label{f:pascal}
\end{figure}

Veremos en la Secci\'on~\ref{s:symmfn} que existen
matrices {\em infinitas} importantes que son totalmente no negativas.
Aunque no est\'a claro que podamos
generalizar el Teorema~\ref{t:tnnconverse} al caso infinito,
podemos generalizar 
el Teorema~\ref{t:wlinlem}. 
Por ejemplo la matriz
$D = [d_{ij}]$ definida como
\begin{equation*}
d_{ij} = \begin{cases} \binom{i}{j}  &\text{si $j \leq i$},\\
                                    0 &\text{si no}
         \end{cases}
\end{equation*}
es totalmente no negativa y es la matriz de caminos de la red plana
infinita de la Figura~\ref{f:pascal}~\cite{GVBin}.  
(Las aristas tienen orientaci\'on hacia la derecha y hacia arriba.)

\section{\textsf{Funciones sim\'etricas}\label{s:symmfn}}

En varias situaciones donde se usan matrices totalmente no negativas
tambi\'en se usan {\em funciones sim\'etricas}.
Un polinomio $f(x) = f(x_1,\dotsc,x_n)$
recibe el nombre de {\em funci\'on sim\'etrica} si satisface que
\begin{equation*}
f(x_1,\dotsc,x_{i-1},x_{i+1},x_i,x_{i+2},\dotsc,x_n)
= f(x)
\end{equation*}
para $i=1,\dotsc,n-1$.
Equivalentemente, el polinomio $f(x)$ es una funci\'on sim\'etrica 
si cualquier reordenamiento de las variables deja el polinomio invariante.
(Para mayor informaci\'on, vea~\cite[Cap.\,1]{M1}, \cite[Cap.\,4]{Sag}, 
\nocite{StanEC1}\cite[Cap.\,7]{StanEC2}.) 
Note que nuestra terminolog\'\i a es un poco enga\~nosa:
una expresi\'on como $\sqrt{x_1x_2}$, que es sim\'etrica pero no es polinomial,
no recibe el nombre de {\em funci\'on sim\'etrica} 
aunque sea una funci\'on con la propiedad de simetr\'\i a.

Unos ejemplos sencillos de funciones sim\'etricas son las
{\em funciones elementales},
\begin{align*}
e_1(x) &= x_1 + x_2 + \cdots = \sum_{i=1}^n x_i,\\
e_2(x) &= x_1x_2 + x_1x_3 + x_2x_3 + \cdots = \sum_{i=1}^{n-1}\sum_{j=i+1}^n x_ix_j,\\
&\vdots \\
e_n(x) &= x_1\cdots x_n,\\
e_{n+1}(x) &= 0,\\
&\vdots
\end{align*}
las {\em funciones homog\'eneas (completas)},
\begin{align*}
h_1(x) &= x_1 + x_2 + \cdots = \sum_{i=1}^n x_i,\\
h_2(x) &= x_1^2 + x_1x_2 + x_2^2 + x_1x_3 + x_2x_3 + x_3^2 + \cdots
= \sum_{i=1}^n\sum_{j=i}^n x_ix_j,\\
&\vdots \\
h_m(x) &= \sum_{1 \leq i_1 \leq \cdots i_m \leq n} x_{i_1} \cdots x_{i_m},\\
&\vdots
\end{align*}
y las 
{\em funciones de sumas de potencias},
\begin{align*}
p_1(x) &= x_1 + x_2 + \cdots = \sum_{i=1}^n x_i,\\
p_2(x) &= x_1^2 + x_2^2 + x_3^2 + \cdots
= \sum_{i=1}^n x_i^2,\\
&\vdots \\
p_m(x) &= \sum_{i=1}^n x_i^m.\\
&\vdots
\end{align*}
Note que el nombre {\em funci\'on homog\'enea} tambi\'en es un poco enga\~noso:
todas las funciones arriba son polinomios homog\'eneos.
A veces conviene utilizar un 
vector infinito  $x= (x_1, x_2,\dotsc)$ de variables; es claro como definir las funciones anteriores en este caso. 
Así tenemos que $e_m(x)\neq 0$ para
todo valor de $m$.

Cualquier funci\'on sim\'etrica $f(x)$ 
se puede expresar como una combinaci\'on 
lineal de productos de funciones elementales,  de funciones homog\'eneas
o de funciones de sumas de potencias.
En cada caso, esta expresi\'on siempre es \'unica.
(Equivalentemente, estas clases de funciones forman tres bases para el
anillo de funciones sim\'etricas). 
Si $f$ tiene coeficientes enteros, entonces 
cuando escribimos $f$
como una combinaci\'on lineal de productos de funciones elementales
o de funciones homog\'eneas, los coeficientes que ob\-te\-nemos son 
n\'umeros enteros.
Por otra parte, si escribimos $f$ como una combinaci\'on lineal de productos
de sumas de potencias, los coeficientes que obtenemos son n\'umeros
racionales.
Por ejemplo tenemos que
\begin{align*}
e_4 &= -h_4 + 2h_3h_1 + h_2^2 -3 h_2h_1^2 + h_1^4, \\
&= -\frac14 p_4 +\frac13 p_3p_1 +\frac18 p_2^2 
-\frac14 p_2p_1^2 +\frac1{24} p_1^4,
\end{align*}
\begin{align*}
h_4 &= -e_4 + 2e_3e_1 + e_2^2 -3 e_2e_1^2 + e_1^4, \\
&= \frac14 p_4 +\frac13 p_3p_1 +\frac18 p_2^2 
+ \frac14 p_2p_1^2 +\frac1{24} p_1^4,
\end{align*}
\begin{align*}
p_4 &= -4e_4 + 4e_3e_1 + 2e_2^2 - 4e_2e_1^2 + e_1^4, \\
&= 4h_4 - 4h_3h_1 - 2h_2^2 + 4h_2h_1^2 - h_1^4.
\end{align*}

\begin{figure}[h]
\centerline{\psfig{figure=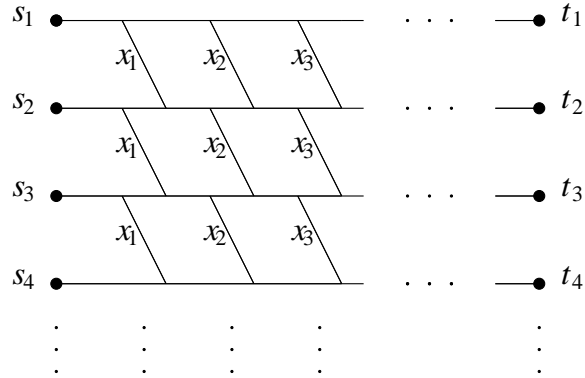,width=.50\textwidth}}
\caption{La red plana para las funciones sim\'etricas elementales.}
\label{f:eplanarnet}
\end{figure}

Cada una de estas tres clases de funciones 
tiene asociada una 
matriz totalmente no negativa.
A la clase de las funciones elementales le asociamos la 
matriz de Toeplitz infinita
\begin{equation}\label{eq:etnn}
E = 
\begin{bmatrix}
1      & e_1    & e_2    & e_3    & e_4    & \ldots \\
0      &   1    & e_1    & e_2    & e_3    & \ldots \\
0      &   0    &   1    & e_1    & e_2    & \ldots \\
0      &   0    &   0    &   1    & e_1    & \ldots \\
0      &   0    &   0    &   0    &  1     & \ldots \\
\vdots &        &        &        & \ddots & \ddots   
\end{bmatrix}.
\end{equation}
(Una matriz $A = [a_{ij}]$ recibe el nombre de {\em Toeplitz} si 
satisface $a_{i,j} = a_{i+k,j+k}$ para cada $k \geq 1$.)
Esta matriz es totalmente no negativa en el sentido de que
cada menor es una funci\'on sim\'etrica con coeficientes positivos.
Una tal funci\'on es denominada {\em monomio-positiva}. 
Si el vector de variables
$x = (x_1,\dotsc,x_n)$ 
es finito, entonces la evaluaci\'on de $E$ en $n$ n\'umeros
no negativos resulta ser una matriz totalmente no negativa
en el sentido corriente.  La matriz $E$ es
la matriz de pesos de la red plana de
la Figura~\ref{f:eplanarnet}.

\begin{figure}[h]
\centerline{\psfig{figure=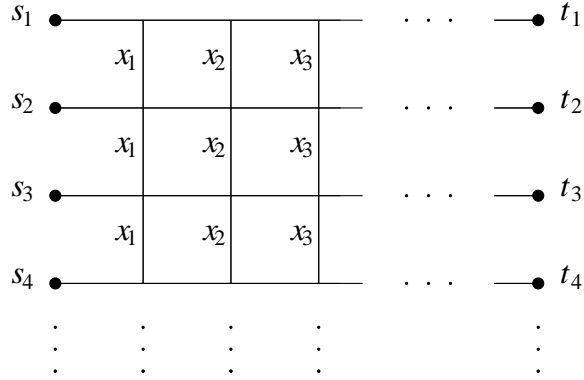,width=.50\textwidth}}
\caption{La red plana para las funciones sim\'etricas homog\'eneas.}
\label{f:hplanarnet}
\end{figure}

A la clase de las funciones homog\'eneas, le asociamos 
la matriz de Toeplitz
\begin{equation}\label{eq:htnn}
H = 
\begin{bmatrix}
1      & h_1    & h_2    & h_3    & h_4    & \ldots \\
0      &   1    & h_1    & h_2    & h_3    & \ldots \\
0      &   0    &   1    & h_1    & h_2    & \ldots \\
0      &   0    &   0    &   1    & h_1    & \ldots \\
0      &   0    &   0    &   0    &  1     & \ldots \\
\vdots &        &        &        & \ddots & \ddots   
\end{bmatrix}.
\end{equation}
Esta matriz tambi\'en es totalmente no negativa, ya que 
es la matriz de pesos de la red plana de la Figura~\ref{f:hplanarnet}, donde las aristas est\'an dirigidas hacia la derecha y hacia abajo.

Para asociar una matriz totalmente no negativa a las funciones de 
sumas de potencias,
necesitamos que el vector de variables $x = (x_1,\dotsc,x_n)$ sea finito.
La matriz de Hankel
\begin{equation}\label{eq:hankel}
P = 
\begin{bmatrix}
n      & p_1    & p_2    & \ldots & p_{n-1} \\
p_1    & p_2    & p_3    & \ldots & p_n     \\
p_2    & p_3    & p_4    & \ldots & p_{n+1} \\
\vdots & \vdots & \vdots &        & \vdots  \\
p_{n-1}& p_n    & p_{n+1}& \ldots & p_{2n-2}    
\end{bmatrix}
\end{equation}
es totalmente no negativa para cualquier conjunto de valores no negativos
$x_1,\dotsc,x_n$.
(Una matriz $A = [a_{ij}]$ recibe el nombre de {\em Hankel} 
si satisface $a_{i,j} = a_{i+k,j-k}$ para $k = 1,\dotsc,j-1$.  
Vea~\cite[Cap.\,10]{Gant1}.)
La matriz $P$ es igual al producto de una matriz de 
Vandermonde con su transpuesta,
$P = VV^T$, donde
\begin{equation}\label{eq:vandermonde}
V = 
\begin{bmatrix}
1         & 1         & \ldots & 1 \\
x_1       & x_2       & \ldots & x_n \\
x_1^2     & x_2^2     & \ldots & x_n^2 \\
\vdots    & \vdots    &        & \vdots \\
x_1^{n-1} & x_2^{n-1} & \ldots & x_n^{n-1}
\end{bmatrix}.
\end{equation}
Si reemplazamos los variables por n\'umeros reales que satisfacen
\begin{equation}\label{eq:increasingpvars}
0 \leq x_1 \leq \cdots \leq x_n,
\end{equation} 
se puede mostrar por inducci\'on que $V$ es totalmente no negativa.
(Vea por ejemplo~\cite[p.\,99]{Gant2}.)  
Concluimos entonces que la matriz $P$ tambi\'en es 
totalmente no negativa.

\begin{figure}[h]
\centerline{\psfig{figure=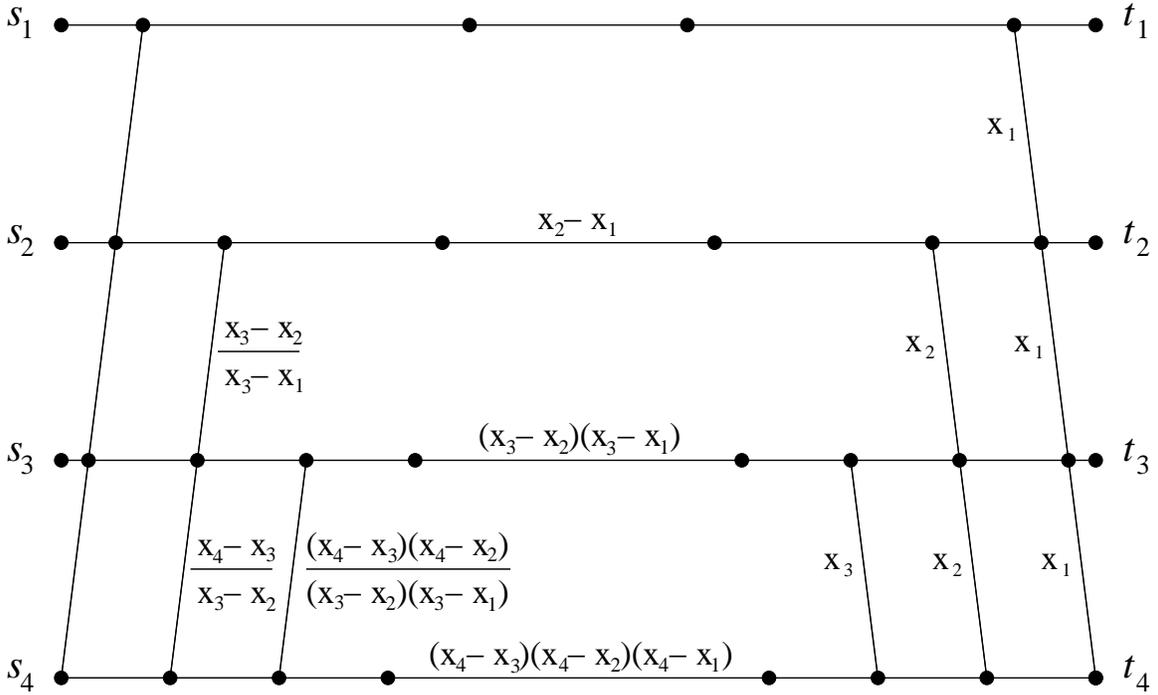,width=\textwidth}}
\caption{Una red plana para la matriz de Vandermonde.}
\label{f:pplanarnet}
\end{figure}

Cuando las desigualdades (\ref{eq:increasingpvars}) son estrictas, 
 la inducci\'on anterior
muestra que $V$ y $P$ son totalmente {\em positivas}.
En este caso es f\'acil factorizar $V$ para construir una red plana.
Ilustraremos esta construcci\'on en el caso $n=4$.
Tenemos que 
\begin{equation*}
V = L_1L_2L_3 D U_3 U_2 U_1, 
\end{equation*}
donde
\begin{equation*}
L_1L_2L_3 =
\begin{bmatrix}
1 & 0 & 0 & 0 \\
1 & 1 & 0 & 0 \\
1 & 1 & 1 & 0 \\
1 & 1 & 1 & 1 
\end{bmatrix}
\begin{bmatrix}
1 & 0 & 0 & 0 \\
0 & 1 & 0 & 0 \\
0 & \tfrac{x_3-x_2}{x_2-x_1} & 1 & 0 \\
0 & \tfrac{x_4-x_3}{x_2-x_1} & \tfrac{x_4-x_3}{x_3-x_2} & 1 
\end{bmatrix}
\begin{bmatrix}
1 & 0 & 0 & 0 \\
0 & 1 & 0 & 0 \\
0 & 0 & 1 & 0 \\
0 & 0 & \tfrac{(x_4-x_3)(x_4-x_2)}{(x_3-x_2)(x_3-x_1)} & 1 
\end{bmatrix},
\end{equation*}
\begin{equation*}
D = 
\begin{bmatrix}
1 & 0 & 0 & 0 \\
0 & x_2-x_1 & 0 & 0 \\
0 & 0 & (x_3-x_2)(x_3-x_1) & 0 \\
0 & 0 & 0 & (x_4-x_3)(x_4-x_2)(x_4-x_1)
\end{bmatrix},
\end{equation*}
\begin{equation*}
U_3U_2U_1 = 
\begin{bmatrix}
1 & 0 & 0 & 0 \\
0 & 1 & 0 & 0 \\
0 & 0 & 1 & x_3 \\
0 & 0 & 0 & 1
\end{bmatrix}
\begin{bmatrix}
1 & 0 & 0 & 0 \\
0 & 1 & x_2 & x_2^2 \\
0 & 0 & 1 & x_2 \\
0 & 0 & 0 & 1
\end{bmatrix}
\begin{bmatrix}
1 & x_1 & x_1^2 & x_1^3 \\
0 & 1 & x_1 & x_1^2 \\
0 & 0 & 1 & x_1 \\
0 & 0 & 0 & 1
\end{bmatrix}.
\end{equation*}
Concatenando las redes planas correspondientes, construimos 
la red plana de la
Figura~\ref{f:pplanarnet}.

Las tres matrices totalmente no negativas (\ref{eq:etnn}), (\ref{eq:htnn}),
(\ref{eq:hankel})
tienen aplicaciones en el problema de localizar las ra\'{\i}ces de un polinomio.
(El n\'umero $\alpha$ es denominado una {\em ra\'{\i}z} o un \emph{cero} del polinomio $a(z)$
si tenemos que $a(\alpha) = 0$.)
 
En los teoremas siguientes, factorizaremos un polinomio real
$a(z) = 1 + a_1 z + \cdots + a_nz^n$ como
\begin{equation}\label{eq:myfactor}
a(z) = \prod_{i=1}^n ( 1 + \beta_i z).
\end{equation}
Los n\'umeros complejos $\beta_1,\dotsc,\beta_n$ 
son los inversos aditivos de los rec\'\i procos de las ra\'{\i}ces de $a(z)$.
Por lo tanto, ellos son reales si y solamente si
los ceros de $a(z)$ son reales.

El resultado siguiente es un caso especial de un teorema de
Aissen, Schoenberg y  
Whitney~\cite{ASW}.
\begin{thm}\label{t:ASW}
Sea $a(z) = 1 + a_1 z + \cdots + a_nz^n$ un polinomio con coeficientes
positivos.  Los ceros de $a(z)$
son reales si y solamente si la matriz de Toeplitz
\begin{equation*}
A = 
\begin{bmatrix}
1      & a_1    & a_2    & a_3    & a_4    & \ldots \\
0      & 1      & a_1    & a_2    & a_3    & \ldots \\
0      &   0    & 1      & a_1    & a_2    & \ldots \\
0      &   0    &   0    & 1      & a_1    & \ldots \\
0      &   0    &   0    & 0      &  1     & \ldots \\
\vdots &        &        &        & \ddots & \ddots   
\end{bmatrix}
\end{equation*}
es totalmente no negativa, donde definimos $a_k=0$ para $k>n$.
\end{thm}
Es inmediato que se puede reemplazar a la matriz $A$
por su transpuesta en el teorema anterior.

Note que cada coeficiente $a_i$ es igual a la funci\'on elemental
$e_i(\beta_1,\dotsc,\beta_n)$.  Por lo tanto la Figura~\ref{f:eplanarnet} (donde cambiamos cada peso $x_i$ por $\beta_i$) 
nos da una red plana cuya matriz de pesos es $A$.
Por otra parte, esta red plana no nos sirve para mostrar que los ceros de 
$a(z)$ son reales, ya que necesitamos saber cu\'ales son los ceros para
poder dibujar la red.

\begin{figure}[h]
\centerline{\psfig{figure=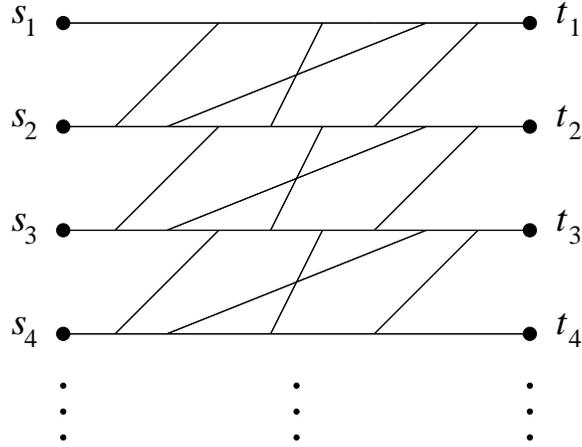,width=.5\textwidth}}
\caption{Una demostraci\'on de que el polinomio $1 + 6z + 5z^2 + z^3$ 
tiene todos los ceros reales.}
\label{f:aswplanarnet}
\end{figure}

Afortunadamente en algunos casos podemos 
construir una red plana para la matriz
$A$ sin saber cu\'ales son los ceros de $a(z)$.
Por ejemplo, la Figura~\ref{f:aswplanarnet} muestra una red plana sin pesos
que demuestra que el polinomio
$1 + 6z + 5z^2 + z^3$ tiene todos los ceros reales
(sin usar los valores de ellos).

La dificultad con este m\'etodo 
es que tenemos que adivinar 
la forma de la red plana.  Adem\'as, no se sabe si existe una tal
red plana sin pesos, ni siquiera en el caso en que
el polinomio tenga todos sus coeficientes enteros.

\begin{question}\label{q:inttoeplitz}
Sea $a(z)$ un polinomio en $\mathbb{N}[z]$ con todos los ceros reales.
?`Cu\'ando y c\'omo es posible
construir una red plana sin pesos para 
mostrar que $a(z)$ tiene todos los ceros reales?
\end{question}

No est\'a claro que sea suficiente multiplicar el polinomio
$a(z)$ por un entero $k$ para poder responder a 
la Pregunta~\ref{q:inttoeplitz} afirmativamente, como sugiere la Observaci\'on~\ref{o:fomin}.
Equivalentemente, no se sabe si existe una red plana
con pesos {\em racionales} para mostrar que $a(z)$ tiene todos los ceros
reales.  El \'unico teorema de factorizaci\'on
de matrices infinitas de Toeplitz que son totalmente no negativas
emplea n\'umeros {\em reales}~\cite{ASW}.

\begin{question}
Sea $a(z)$ un polinomio en $\mathbb{N}[z]$ con todos los ceros reales.
?`Cu\'ando y c\'omo es posible
construir una red plana con pesos racionales para 
mostrar que $a(z)$ tiene todos los ceros reales?
\end{question}

De manera similar al Teorema~\ref{t:ASW}, el siguiente teorema
emplea a las funciones de sumas de potencias para demostrar
que un polinomio tiene todos sus ceros reales.
Este teorema es una consecuencia 
de un 
resultado de Gantmacher~\cite[Cor.\,de Tm.\,6,~p.\,203]{Gant2}, 
en combinaci\'on con
la positividad total de la matriz de Vandermonde~(\ref{eq:vandermonde}).  
(Vea tambi\'en~[Tm.\,6.5]\cite{Skanhsurvey}.)
\begin{thm}\label{t:gant}
Sea $a(z) = 1 + a_1 z + \cdots + a_nz^n$
un polinomio con coeficientes positivos.
Todos los ceros de $a(z)$
son reales y distintos si y solamente si la matriz de Hankel $P$ 
de~(\ref{eq:hankel}) es totalmente positiva, donde
$p_i$ es la suma de potencias $\beta_1^i + \cdots + \beta_n^i$,
y $\beta_1,\dotsc,\beta_n$ son los n\'umeros que aparecen en la
factorizaci\'on (\ref{eq:myfactor}).
\end{thm}
Sin saber cu\'ales son los valores de $\beta_1,\dotsc,\beta_n$,
es posible calcular la sumas de potencias en t\'erminos de los
coeficientes $a_1,\dotsc,a_n$, que como ya hemos observado son
las funciones 
sim\'etricas elementales evaluadas en $\beta_1,\dotsc,\beta_n$.
Otra vez es concebible que podamos construir una red plana
sin pesos para mostrar que $a(z)$ tiene todos los ceros reales.
Desa\-fortunadamente no se conoce un m\'etodo general para realizar
esta construcci\'on, aun cuando la matriz $P$ tiene entradas enteras
y la forma especial de Hankel.

\begin{question}
Sea $A$ una matriz totalmente no negativa con entradas enteras.
?`Cu\'ando y c\'omo es posible realizar $A$ como la matriz de
caminos de una red plana sin pesos?
?`C\'omo se puede encontrar el menor entero $k$ de modo que
$kA$ sea realizable as\'\i ?
?`Es m\'as f\'acil el caso especial de matrices de Hankel?
\end{question}

El resultado de Gantmacher~\cite[Cor.\,de Tm.\,6, p.\,203]{Gant2} 
se puede modificar para aplicarse a los polinomios
con ceros repetidos~\cite[Cap.\,15, Sec.\,11]{Gant2}.  
Ser\'\i a interesante modificar el Teorema~\ref{t:gant}
tambi\'en.

\begin{question}
Sea
$a(z) = 1 + a_1 z + \cdots + a_nz^n$ un polinomio con coeficientes
positivos.
?`Es cierto que todos los ceros de $a(z)$
son reales
(y no necesariamente distintos) si y solamente si
la matriz de Hankel $P$ (\ref{eq:hankel}) es totalmente no negativa?
\end{question}

\section{\textsf{Polinomios totalmente no negativos}\label{s:ptnn}}

Podemos considerar al determinante de una matriz $n \times n$ como un polinomio
en $n^2$ variables $x = ( x_{11}, \dotsc, x_{nn})$
\begin{equation*}
\det(x) = \sum_{\sigma \in S_n} \text{sgn}( \sigma ) 
x_{1,\sigma(1)} \cdots x_{n,\sigma(n)}
\end{equation*}
que se puede evaluar en una matriz $A = [a_{ij}]$
usando la sustituci\'on $x_{ij} = a_{ij}$.
Seg\'un la definici\'on de las matrices totalmente no negativas,
la evaluaci\'on de este polinomio en una matriz totalmente no negativa
resulta siempre en un n\'umero no negativo.
Un polinomio con esta propiedad se denomina
{\em totalmente no negativo}.
Los polinomios totalmente no negativos aparecen en la investigaci\'on de
Lusztig~\cite{LusztigTP}, quien mostr\'o que todos los elementos de la
base can\'onica del anillo de coordenadas de $GL_n$ son polinomios totalmente
no negativos.
Para comprender mejor esta base, que todav\'\i a no tiene una descripci\'on
sencilla, quiz\'as uno podr\'\i a esperar 
caracterizar todos los polinomios
totalmente no negativos, o por lo menos un subconjunto de ellos.

\begin{figure}[h]
\centerline{\psfig{figure=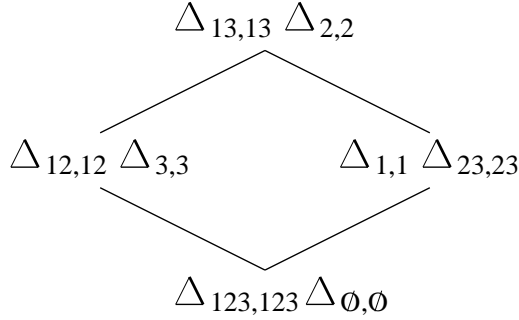,width=.45\textwidth}}
\caption{Un \'orden parcial sobre los productos de menores principales
de matrices 
$3 \times 3$ 
totalmente no negativas.}
\label{f:tnnposet3}
\end{figure}

Fallat y sus colaboradores descubrieron un subconjunto de los polinomios
totalmente no negativos que 
se puede describir en t\'erminos de menores principales.
Un ejemplo es el polinomio 
\begin{equation}\label{eq:1324diff}
\Delta_{\{1,3\}\{1,3\}}
\Delta_{\{2\}\{2\}}
-
\Delta_{\{1\}\{1\}}
\Delta_{\{2,3\}\{2,3\}}.
\end{equation}
Equivalentemente todas las matrices totalmente no negativas de tama\~no
por lo menos $3 \times 3$ satisfacen la desigualdad
\begin{equation}\label{eq:13,24}
\Delta_{\{1\}\{1\}}
\Delta_{\{2,3\}\{2,3\}}
\leq
\Delta_{\{1,3\}\{1,3\}}
\Delta_{\{2\}\{2\}}.
\end{equation}
Esta desigualdad y cuatro parecidas se muestran 
en la Figura~\ref{f:tnnposet3} en la forma de
un conjunto parcialmente ordenado. 
(Interpretamos el menor vac\'\i o 
$\Delta_{\emptyset,\emptyset}$ como $1$.) 

Los Teoremas~\ref{t:linlem} y \ref{t:tnnconverse}  
proveen m\'etodos para demostrar que los productos
\begin{equation}\label{eq:products}
\Delta_{13,13}\Delta_{2,2},\quad
\Delta_{12,12}\Delta_{3,3},\quad
\Delta_{1,1}\Delta_{23,23},\quad
\Delta_{123,123}\Delta_{\emptyset, \emptyset}
\end{equation}
de menores de matrices totalmente no negativas
est\'an relacionados como indica la Figura~\ref{f:tnnposet3}.
Para demostrarlo, observemos primero que estos teoremas implican el siguiente corolario.

\begin{cor}\label{c:prod}
Sea $A$ una matriz $n \times n$ totalmente no negativa,
y sea $G$ una red plana cuya matriz de caminos es $A$.
Sea $I$ un subconjunto de $[n]$ y sea $\notI = [n] \ssm I$
su complemento en $[n]$.
Entonces el producto $\Delta_{I,I}\Delta_{\notI,\notI}$ 
de menores de $A$ es igual a la suma de pesos de familias
$\pi = (\pi_1,\dotsc,\pi_n)$ de caminos en $G$ con las 
siguientes
caracter\'\i sticas.
\begin{enumerate}
\item En cada camino, el \'\i ndice de la fuente es igual al
\'\i ndice del destino.
\item Los caminos con \'\i ndices en $I$ no se intersectan entre ellos,
y los caminos con \'\i ndices en $\notI$ no se intersectan entre
ellos.
\end{enumerate}
\end{cor}

Imaginando que las fuentes y los destinos con \'\i ndices
$I$ son de un color y que las fuentes y los destinos con \'\i ndices
$\notI$ son de otro color, podemos interpretar 
el producto $\Delta_{I,I}\Delta_{\notI,\notI}$ 
como el n\'umero de familias de caminos
en las cuales cada camino tiene uno de dos colores
y dos caminos del mismo color no se pueden intersectar.

\begin{figure}[h]
\centering
\mbox{\subfigure[Contado en $\Delta_{123,123}\Delta_{\emptyset,\emptyset}$.]
     {\epsfig{figure=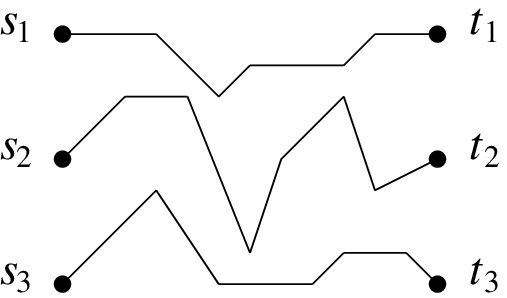,width=.35\textwidth}
      \label{f:pretlieb3idfams:a}} \qquad\qquad
      \subfigure[Contado en $\Delta_{13,13}\Delta_{2,2}$.]
     {\epsfig{figure=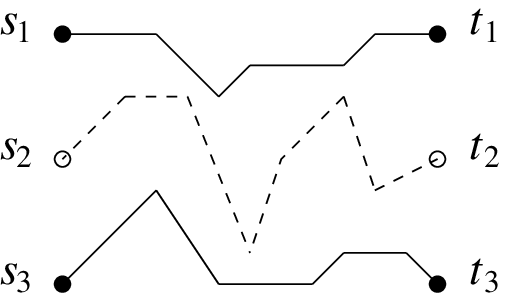,width=.35\textwidth}
      \label{f:pretlieb3id1}} } \\
\mbox{\subfigure[Contado en $\Delta_{12,12}\Delta_{3,3}$.]
     {\psfig{figure=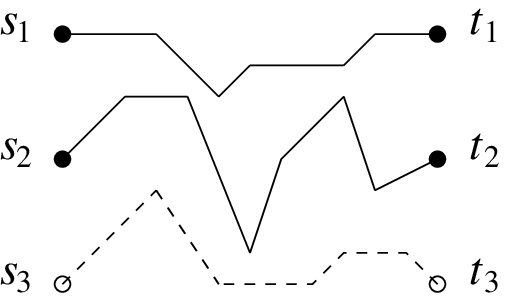,width=.35\textwidth}
      \label{f:pretlieb3id2}} \qquad\qquad
      \subfigure[Contado en $\Delta_{1,1}\Delta_{23,23}$.]
     {\psfig{figure=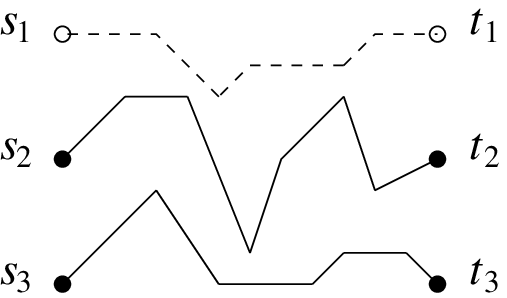,width=.35\textwidth}
      \label{f:pretlieb3id3}} }
\caption{Coloraciones de una familia de tres caminos.}
\label{f:pretlieb3idfams}
\end{figure}

Ahora consideremos los conjuntos de aristas que son contados
por un producto
$\Delta_{I,I}\Delta_{\notI,\notI}$ pero no por otro.
La Figura~\ref{f:pretlieb3idfams}(a) muestra un conjunto de aristas
que define una familia de caminos sin intersecci\'on.
Este conjunto de aristas es contado por todos los cuatro productos
(\ref{eq:products}), como se ve facilmente en las 
Figuras~\ref{f:pretlieb3idfams}(a),
\ref{f:pretlieb3idfams}(b),
\ref{f:pretlieb3idfams}(c), y
\ref{f:pretlieb3idfams}(d).

\begin{figure}[h]
\centering
\mbox{\subfigure[Un conjunto de aristas que puede definir tres caminos]
     {\psfig{figure=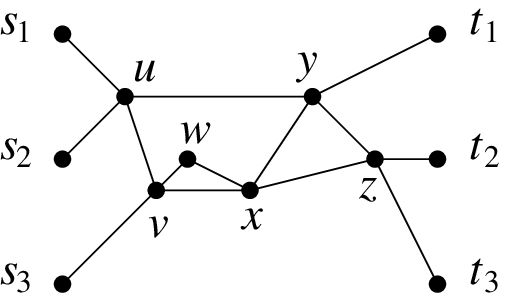,width=.35\textwidth}
      \label{f:pretlieb3al}} \qquad\qquad
      \subfigure[Las restricciones de los colores en una coloraci\'on del grafo anterior.]
     {\psfig{figure=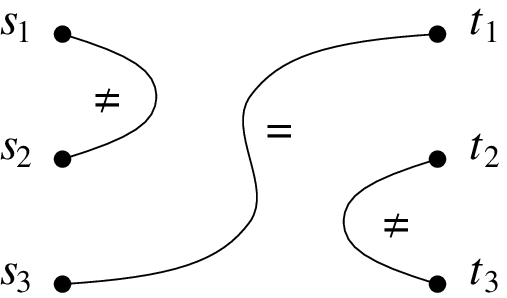,width=.35\textwidth}
      \label{f:tlieb3al}} } \\
\mbox{\subfigure[Una familia coloreada de tres caminos]
     {\psfig{figure=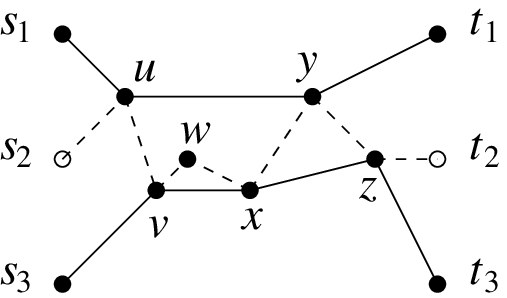,width=.35\textwidth}
      \label{f:pretlieb3al1}} \qquad\qquad
      \subfigure[Otra familia coloreada de tres caminos]
     {\psfig{figure=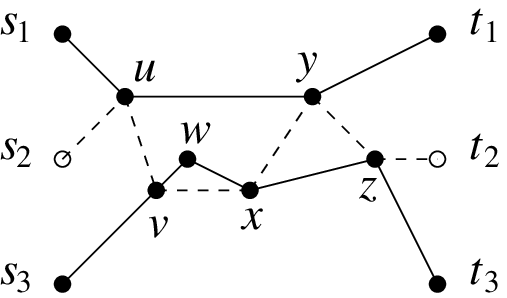,width=.35\textwidth}
      \label{f:pretlieb3al2}} } \\
\caption{Familias de caminos.}
\label{f:pretlieb3alfams}
\end{figure}

Por otra parte, el conjunto de aristas en la Figura~\ref{f:pretlieb3alfams}(a)
no es contado por todos los cuatro productos~(\ref{eq:products}).
Para que este conjunto de aristas determine una familia $2$-coloreable
de caminos, las aristas $(s_1,u)$ y 
$(s_2,u)$ necesitan tener colores distintos, porque se intersectan en $u$.
Entonces $s_1$ y $s_2$ tienen colores distintos.
Por la misma raz\'on $t_2$ y $t_3$ tienen colores distintos.
Mirando las dem\'as aristas, vemos que las aristas en la sucesi\'on
\begin{equation*}
(s_3,v), (u,v), (u,y), (x,y), (x,z), (y,z), (z,t_1)
\end{equation*}
necesitan alternar de color.
Por lo tanto $s_3$ y $t_1$ deben tener el mismo color.
La Figura~\ref{f:pretlieb3alfams}(b) muestra 
esta igualdad y las dos desigualdades.

No es dif\'\i cil mostrar que, al llevar a cabo este an\'alisis para cualquier grafo  que corresponde a una familia $2$-coloreable de tres caminos, el resultado pertenece siempre a una
de las cinco categor\'\i as
siguientes~\cite[Sec.\,2]{SkanIneq}:

%
\begin{equation*}
\psfig{figure=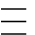,width=1cm}, \quad
\psfig{figure=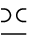,width=1cm}, \quad
\psfig{figure=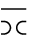,width=1cm}, \quad
\psfig{figure=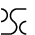,width=1cm}, \quad
\psfig{figure=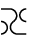,width=1cm}.
\end{equation*}
Estos cinco diagramas
se\~nalan parejas de fuentes o de destinos que
necesitan tener distintos colores, y parejas de una fuente y un destino
que necesitan tener el mismo color.

(Vale la pena anotar que los cinco diagramas son una base
del \'algebra de {\em Temperley-Lieb} $T_3$, 
si definimos el producto de dos diagramas como su concatenaci\'on.
En general formamos el \'algebra de Temperley-Lieb $T_n$
con los $\tfrac{1}{n+1}\tbinom{2n}{n}$ diagramas an\'alogos en $2n$ v\'ertices. Llamamos a estos diagramas la \emph{base est\'andar} de $T_n$.
El {\em n\'umero de Catal\'an}
$\tfrac{1}{n+1}\tbinom{2n}{n}$ aparece con gran frecuencia en la combinatoria algebraica.
Vea por ejemplo la Secci\'on~\ref{s:fed2} y \cite[pp.\,219-231]{StanEC2}.)


Entonces, para comparar productos de la forma
$\Delta_{I,I}\Delta_{\notI,\notI}$ 
en matrices totalmente no negativas
$3 \times 3$,
podemos usar subconjuntos de la base est\'andar del \'algebra de Temperley-Lieb $T_3$.
Para cada producto
$\Delta_{I,I}\Delta_{\notI,\notI}$, hacemos lo siguiente.
\begin{enumerate}
\item Dibujamos dos columnas de $n$ v\'ertices. 
Llamamos a los v\'ertices a la izquierda {\em fuentes} 
y a los v\'ertices a la derecha {\em destinos}.  
(Los denotamos $s_1,\dotsc,s_n$ y $t_1,\dotsc,t_n$ como siempre.)
\item Para $i=1,\dotsc,n$, coloreamos $s_i$ y $t_i$ de rojo
si $i \in I$,  y de azul si $i \in \notI$.
\item Apuntamos los elementos del \'algebra de Temperley-Lieb $T_n$
que conectan solamente las parejas siguentes.
\begin{enumerate}
\item dos fuentes de colores distintos.
\item dos destinos de colores distintos.
\item una fuente y un destino del mismo color.
\end{enumerate}
\end{enumerate}

\begin{figure}[h]
\centerline{\psfig{figure=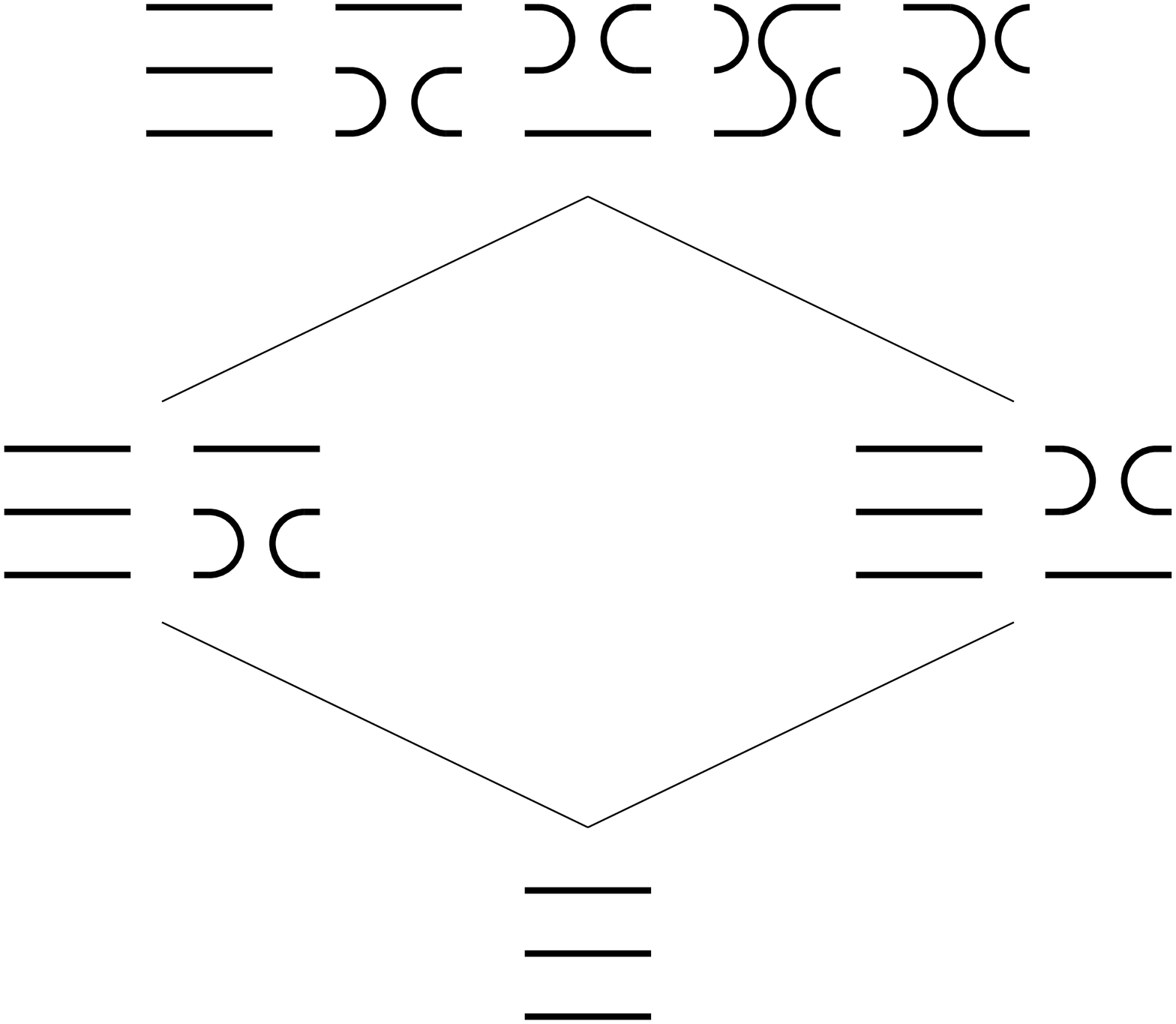,width=.55\textwidth}}
\caption{Un poset de ciertos subconjuntos de la base de $T_3$.}
\label{f:tnntlposet3}
\end{figure}

As\'{\i}, para cada producto $\Delta_{I,I}\Delta_{\notI,\notI}$ 
obtenemos un subconjunto de la base est\'andar de $T_n$.
Ordenando estos subconjuntos por inclusi\'on,
obtenemos todas las desigualdades de la forma
\begin{equation}\label{eq:tnnineq}
\Delta_{I,I}\Delta_{\notI,\notI} \leq 
\Delta_{J,J}\Delta_{\notJ,\notJ}.
\end{equation}
que son v\'alidas para toda matriz totalmente no negativa.
La Figura~\ref{f:tnntlposet3} muestra otra vez
el poset de la Figura~\ref{f:tnnposet3},
ahora con productos de menores reemplazados por subconjuntos de la base est\'andar  de $T_3$.

\begin{figure}
\centerline{\psfig{figure=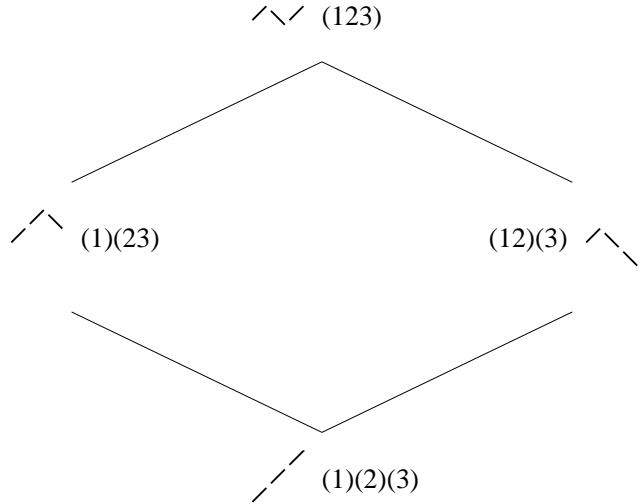,width=.55\textwidth}}
\caption{Un poset de ciertos caminos de tres pasos.}
\label{f:tnnlatticeposet3}
\end{figure}

Hay un criterio m\'as sencillo
para comparar productos de menores de matrices totalmente no 
negativas~\cite{RSkanTLImm}.
\begin{enumerate}
\item Dibujamos un camino de $n$ pasos de longitud uno, de modo que
el paso $i$ tiene pendiente $1$ si $i \in I$ y pendiente $-1$ si $i \in \notI$.
\item Definimos una partici\'on del conjunto $[n]$ de modo que cada bloque
consista de todos los pasos que est\'en a la misma altura.
\end{enumerate}
Ordenando estas particiones de $[n]$ por refinamiento, obtenemos
otra vez todas las desigualdades de la forma (\ref{eq:tnnineq}).
La Figura~\ref{f:tnnlatticeposet3} muestra otra vez el poset 
de la Figura~\ref{f:tnnposet3},
ahora con productos de menores reemplazados por caminos.

Hasta ahora hemos caracterizado todos los polinomios totalmente no negativos
de la forma
\begin{equation*}
\Delta_{J,J}\Delta_{\notJ,\notJ} -
\Delta_{I,I}\Delta_{\notI,\notI}  
\end{equation*}
y estos resultados se pueden generalizar
para obtener una caracterizacion~\cite{FGJMult}, \cite{SkanIneq} 
de todos los polinomios totalmente no negativos de la forma
\begin{equation*}
\Delta_{J,J'}\Delta_{L,L'} -
\Delta_{I,I'}\Delta_{K,K'}.  
\end{equation*}

Ser\'\i a interesante tener una caracterizaci\'on de polinomios
totalmente no negativos m\'as generales.

\begin{question}
?`C\'omo podemos describir combinatoriamente el poset de productos de tres o m\'as menores de matrices totalmente no negativas?
\end{question}

\begin{question}
?`Existe una caracterizaci\'on sencilla de todos los polinomios totalmente
no negativos?
\end{question}

\section{\textsf{Funciones de Schur}\label{s:schur}}

En la Secci\'on~\ref{s:symmfn} hemos considerado
la matriz totalmente no negativa
\begin{equation*}
H = 
\begin{bmatrix}
1      & h_1    & h_2    & h_3    & h_4    & \ldots \\
0      &   1    & h_1    & h_2    & h_3    & \ldots \\
0      &   0    &   1    & h_1    & h_2    & \ldots \\
0      &   0    &   0    &   1    & h_1    & \ldots \\
0      &   0    &   0    &   0    & 1      & \ldots \\
\vdots &        &        &        & \ddots & \ddots 
\end{bmatrix}.
\end{equation*}
Esta matriz, sus submatrices 
y sus menores son importantes en la teor\'\i a de las 
representaciones y en la geometr\'\i a algebraica.
La matriz $H$ se conoce como la {\em matriz infinita de Jacobi-Trudi},  
y sus submatrices tambi\'en se llaman {\em matrices de Jacobi-Trudi}.
Cada menor de $H$ es una funci\'on sim\'etrica,
ya que el conjunto de funciones sim\'etricas es un anillo.
Los menores de $H$ que corresponden a conjuntos {\em  consecutivos} de columnas
se denominan {\em funciones de Schur}. En la segunda parte de esta serie de art\'{\i}culos \cite{ALRS2} tendremos m\'as que decir sobre ellas.

Note que una funci\'on de Schur no tiene un nombre \'unico seg\'un
nuestra notaci\'on $\Delta_{I,I'}$.  Por ejemplo tenemos que
\begin{equation*}
\det
\begin{bmatrix}
h_2 & h_3 & h_4 \\
h_1   & h_2 & h_3 \\
0   & 1   & h_1
\end{bmatrix}
= \Delta_{124,345}
= \Delta_{235,456}
= \Delta_{346,567}
= \cdots
\end{equation*}

Por costumbre, denotamos a cada funci\'on de Schur
usando los \'\i ndices en la diagonal de la submatriz correspondiente de $H$. 
Por ejemplo, la funci\'on de Schur del ejemplo anterior es $s_{221}$. Para cada pareja $I,I'$ tal que $s_{221}=\Delta_{I,I'}$ tenemos que $I'-I = (2,2,1)$ si tratamos a $I$ e $I'$ como vectores. Los \'\i ndices de una funci\'on de Schur siempre forman una sucesi\'on
decreciente (no necesariamente estrictamente), que tambi\'en se llama
una {\em partici\'on}.  Una partici\'on se suele
escribir como
\begin{equation*}
\lambda = (\lambda_1,\dotsc,\lambda_r).
\end{equation*}
Los n\'umeros $\lambda_1,\dotsc,\lambda_r$ se llaman las {\em partes} de
$\lambda$.
Decimos que $\lambda$ es una {\em partici\'on de $n$} si tenemos que
\begin{equation*}
\lambda_1 + \cdots + \lambda_r = n.
\end{equation*}
En este caso tambi\'en escribimos $\lambda \vdash n$, o $|\lambda| = n$.

Recordando que $H$ es la matriz de pesos de la red plana en la
Figura~\ref{f:hplanarnet}, 
podemos interpretar cada funci\'on de Schur como una suma de pesos  
de caminos en ese grafo.
Entonces est\'a claro que la expansi\'on de una funci\'on de Schur en
t\'erminos de las variables $x_1,\dotsc$ no tiene resta.

\begin{obs}
Cada funci\'on de Schur es monomio-positiva.
\end{obs}

\begin{figure}
\centerline{\psfig{figure=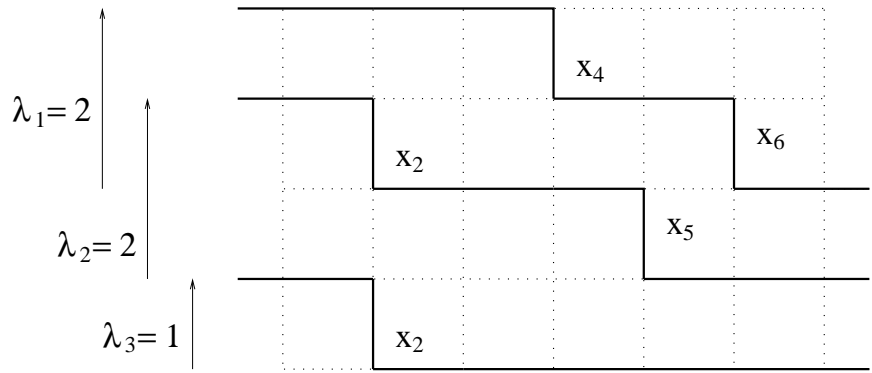,width=.75\textwidth}}
\caption{Una familia de caminos que contribuye $x_2^2x_4x_5x_6$ a 
la funci\'on $s_{221}$.} 
\label{f:splanarnet}
\end{figure}

La Figura~\ref{f:splanarnet} muestra una familia de caminos
que contribuye con un peso de $x_2^2x_4x_5x_6$ a la funci\'on $s_{221}$.

En analog\'\i a con las funciones elementales, homog\'eneas, 
y de sumas de potencias, las funciones
de Schur forman una base del anillo de funciones sim\'etricas, visto como
un espacio vectorial.
Es decir, cualquier funci\'on sim\'etrica $f(x)$ se puede expresar
de manera \'unica como una combinaci\'on lineal de funciones de Schur.
De hecho, si $f$ tiene coeficientes enteros, la combinaci\'on lineal
tambi\'en tendr\'a coeficientes enteros.
Por ejemplo tenemos que
\begin{equation}\label{eq:tos}
\begin{split}
e_4 &= s_{1111}, \\
h_4 &= s_4, \\
p_4 &= s_4 - s_{31} + s_{22} - s_{211}.
\end{split}
\end{equation}
Notemos que aunque la funci\'on $p_4$ es monomio-positiva, no es
Schur-positiva.  Es decir, su expansi\'on en t\'erminos de funciones
de Schur no tiene coeficientes positivos. 

Una propiedad interesante de las funciones de Schur es que el producto
de dos de ellas siempre es Schur-positivo.  La expansi\'on (\'unica) de un
producto tal en t\'erminos de funciones de Schur es 
una combinaci\'on lineal {\em no negativa} de ellas.  Dadas dos particiones
$\lambda$, $\mu$, 
podemos escribir
\begin{equation}\label{eq:lr}
s_\lambda s_\mu = \sum_{\nu \vdash |\lambda| + |\mu|} c_{\lambda\mu}^\nu s_\nu.
\end{equation}
Los coeficientes $c_{\lambda\mu}^\nu$ se llaman los
{\em coeficientes de Littlewood-Richardson}.
Por ejemplo tenemos que
\begin{equation}\label{eq:s31s21}
s_{31}s_{21} = s_{52} + s_{511} + s_{43} + 2s_{421} + s_{4111} + 
s_{331} + s_{322} + s_{3211}.
\end{equation}
Para obtener esta expresi\'on, 
podr\'\i amos en principio
calcular los determinantes $s_{31}, s_{21}$, 
multiplicarlos,
y reconocer la expresi\'on que resulta como una combinaci\'on lineal
de varios otros determinantes.
Afortunadamente, existen m\'etodos menos fastidiosos y m\'as efectivos
de calcular la expansi\'on de un producto $s_\lambda s_\mu$ en t\'erminos de funciones de Schur. 
Todos estos m\'etodos suelen llamarse
{\em el m\'etodo de Littlewood-Richardson}, 
aunque sean diferentes.

En la presentaci\'on de un m\'etodo de Littlewood-Richardson,
vamos a representar a
$\lambda$ y $\mu$ como diagramas de cuadrados que se 
llaman {\em diagramas de Young}. 
Rellenamos los cuadrados del diagrama de $\mu$ con n\'umeros, creando
una tabla que se llama un {\em tableau de Young}.  
Un tableau de Young se llama {\em semiest\'andar} 
si los n\'umeros crecen (no necesariamente estrictamente)
en las filas y crecen estrictamente en las columnas.
Un tableau semiest\'andar de Young de forma $\mu$ se llama
{\em est\'andar} si los n\'umeros usados son
$1,\dotsc,|\mu|$.

Definimos el {\em contenido} de un tableau de Young $T$ como el vector
\begin{equation*}
c(T) = ( \# \text{ de } 1\text{s} \text{ en $T$}, \# \text{ de } 2\text{s} \text{ en $T$},\dotsc)
\end{equation*}

Por ejemplo, un tableau semiest\'andar de Young 
de forma $(4,3,1,1)$ y contenido $(2,2,3,1,1)$
es el siguiente,
\begin{equation}\label{eq:ssyt}
\NantelYoung{\Carre{1}&\Carre{1}&\Carre{2}&\Carre{3}\cr
           \Carre{2}&\Carre{3}&\Carre{3}\cr
           \Carre{4}\cr
           \Carre{5}\cr
           }.
\end{equation}


Si $T$ es un tableau de Young de forma $\mu$, 
escribiremos $T_{j}$ para denotar el tableau que consiste solamente de las
columnas $j$ hasta $\mu_1$ de $T$.  Entonces si $T$ es el tableau 
semiest\'andar de Young en (\ref{eq:ssyt}), tenemos que
\begin{align*}
c(T_1) &= (2,2,3,1,1),\\
c(T_2) &= (1,1,3,0,0),\\
c(T_3) &= (0,1,2,0,0),\\
c(T_4) &= (0,0,1,0,0).
\end{align*}

Uno de los m\'etodos de 
Littlewood-Richardson para multiplicar $s_\lambda s_\mu$ es 
el siguiente.
\begin{enumerate}
\item Dibujamos $\lambda$ y $\mu$ como diagramas de Young.
\item Escribimos todos los tableaux semiest\'andar de Young $T$ 
de forma $\mu$ tales que para $j=1,\dotsc,\mu_1$, se cumpla que
$\lambda + c(T_j)$ es una partici\'on.
\item 
Sumamos $s_{\lambda + c(T)}$ sobre todos los 
tableaux $T$ en la lista.
\end{enumerate}
Veamos por ejemplo la multiplicaci\'on de $s_{31} s_{21}$.
Se puede verificar que los tableaux 
semiest\'andar de Young de la forma $(2,1)$ para los cuales
$(3,1) + c(T)$ y $(3,1) + c(T_2)$ son particiones
son los siguientes:
\begin{equation*}
\NantelYoung{\Carre{1}&\Carre{1}\cr
           \Carre{2}\cr
           },\quad
\NantelYoung{\Carre{1}&\Carre{1}\cr
           \Carre{3}\cr
           },\quad
\NantelYoung{\Carre{1}&\Carre{2}\cr
           \Carre{2}\cr
           },\quad
\NantelYoung{\Carre{1}&\Carre{1}\cr
           \Carre{3}\cr
           },\quad
\NantelYoung{\Carre{1}&\Carre{3}\cr
           \Carre{2}\cr
           },\quad
\NantelYoung{\Carre{1}&\Carre{1}\cr
           \Carre{4}\cr
           },\quad
\NantelYoung{\Carre{2}&\Carre{2}\cr
           \Carre{3}\cr
           },\quad
\NantelYoung{\Carre{2}&\Carre{3}\cr
           \Carre{3}\cr
           },\quad
\NantelYoung{\Carre{2}&\Carre{3}\cr
           \Carre{4}\cr
           }.
\end{equation*}
Entonces las particiones $(3,1) + c(T)$ son
\begin{gather*}
(5,2) \quad (5,1,1) \quad (4,3) \quad (4,2,1) \quad (4,2,1) \\
(4,1,1,1) \quad (3,3,1) \quad (3,2,2) \quad (3,2,1,1).
\end{gather*}
As\'\i~obtenemos la expresi\'on de la Ecuaci\'on~(\ref{eq:s31s21}).

Una generalizaci\'on de las funciones de Schur son los menores
de la matriz $H$ que corresponden a conjuntos de columnas
{\em no} necesariamente consecutivas. 
Estas funciones se llaman {\em funciones sezgadas de Schur} 
y se describen usando una {\em pareja} de particiones: $s_{\lambda/\mu}$.  
La primera partici\'on 
determina un menor $k \times k$ de columnas consecutivas, 
y la segunda partici\'on mueve la primera columna 
$\mu_1$ posiciones a la izquierda,
la segunda columna $\mu_2$ posiciones a la izquierda, etc.
Por ejemplo podemos ver en la submatriz 
\begin{equation*}
\begin{bmatrix}
h_2 & h_3 & h_4 & h_5 & h_6 \\
h_1 & h_2 & h_3 & h_4 & h_5 \\
0   & 0   & 0   & 1   & h_1
\end{bmatrix},
\end{equation*}
que las funciones $s_{441}$ y $s_{441/21}$ son
\begin{equation*}
s_{441} = \det
\begin{bmatrix}
h_4 & h_5 & h_6 \\
h_3 & h_4 & h_5 \\
0   & 1   & h_1
\end{bmatrix},\qquad
s_{441/21} = \det
\begin{bmatrix}
h_2 & h_4 & h_6 \\
h_1 & h_3 & h_5 \\
0   & 0   & h_1
\end{bmatrix}.
\end{equation*}

\begin{figure}
\centerline{\psfig{figure=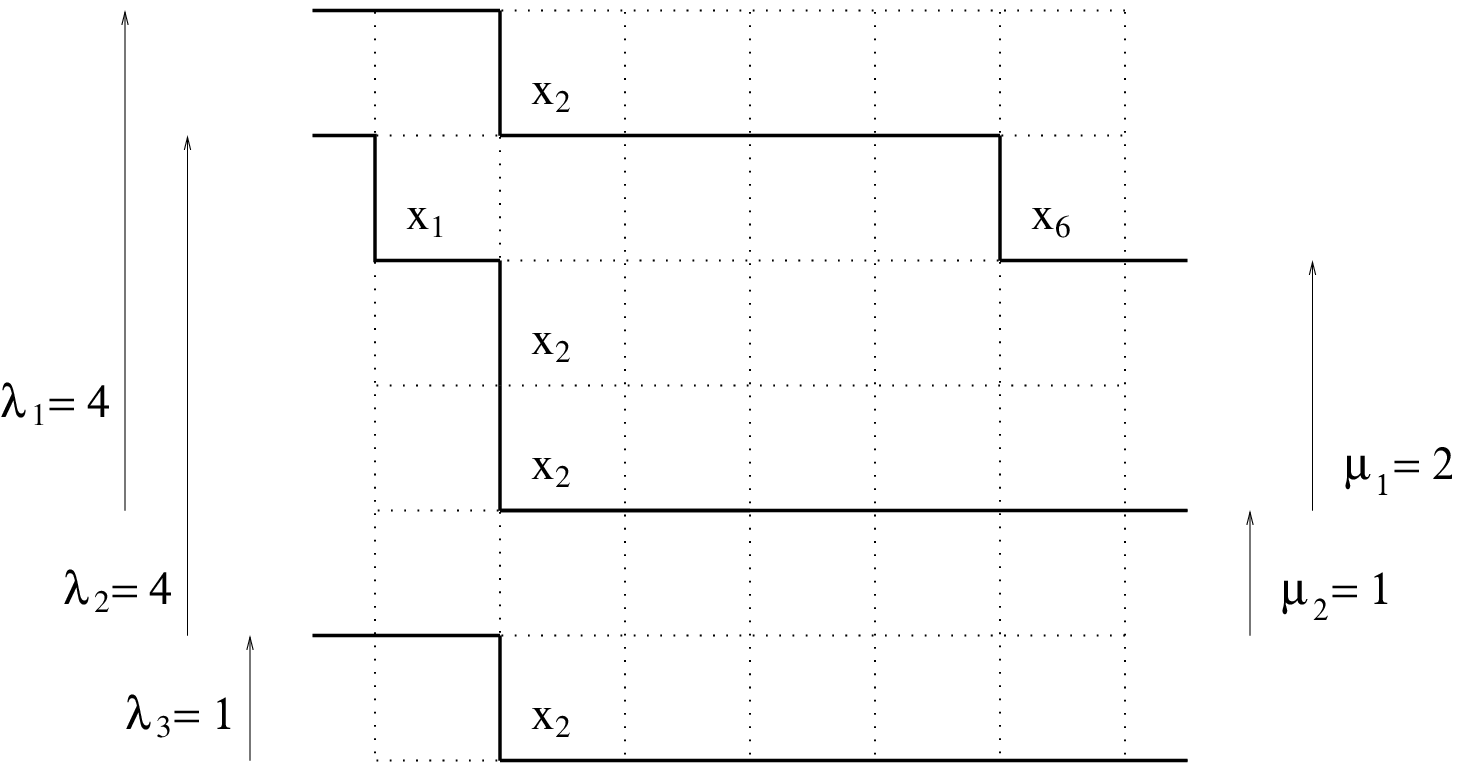,width=.80\textwidth}}
\caption{Una familia de caminos que contribuye 
$x_1x_2^4x_6$ a la funci\'on sezgada de Schur $s_{441/21}$.} 
\label{f:ssplanarnet}
\end{figure}

La Figura~\ref{f:ssplanarnet} muestra
una familia de caminos que contribuye con un peso de $x_1x_2^4x_6$ 
a la funci\'on sezgada de Schur $s_{441/21}$.

Entonces est\'a claro que la expansi\'on de una funci\'on sezgada de Schur en
t\'erminos de las variables $x_1,\dotsc$ no tiene resta.

\begin{obs}
Cada funci\'on sezgada de Schur es monomio-positiva.
\end{obs}

Una funci\'on sezgada de Schur no es solamente monomio-positiva, sino
tambi\'en Schur-positiva.  Además, los coeficientes positivos que aparecen
en la expansi\'on de Schur de una funci\'on sezgada de Schur son los
mismos coeficientes de Littlewood-Richardson que obtenemos en la 
multiplicaci\'on de funciones corrientes de Schur~(\ref{eq:lr}):
\begin{equation}\label{eq:jdt}
s_{\nu/\mu} = \sum_{\lambda \vdash |\nu| - |\mu|}c_{\lambda\mu}^\nu s_\lambda.
\end{equation}

An\'alogo al m\'etodo de Littlewood-Richardson para multiplicar
funciones de Schur,
es el {\em juego de taqu\'\i n} 
que construye la expresi\'on (\ref{eq:jdt}) dadas dos particiones
$\mu$, $\nu$.  Jugamos el juego de taqu\'\i n en un tableau est\'andar 
de Young de la forma $\nu/\mu$.  Para obtener tal tableau,
dibujamos el diagrama de Young para $\nu$, borramos las cajas
que corresponden a $\mu$,
y rellenamos los cuadrados que quedan con los n\'umeros 
$1,\dotsc,|\nu|-|\mu|$ de modo que crezcan en la filas y en las columnas.

\begin{figure}[h]
\centerline{\psfig{figure=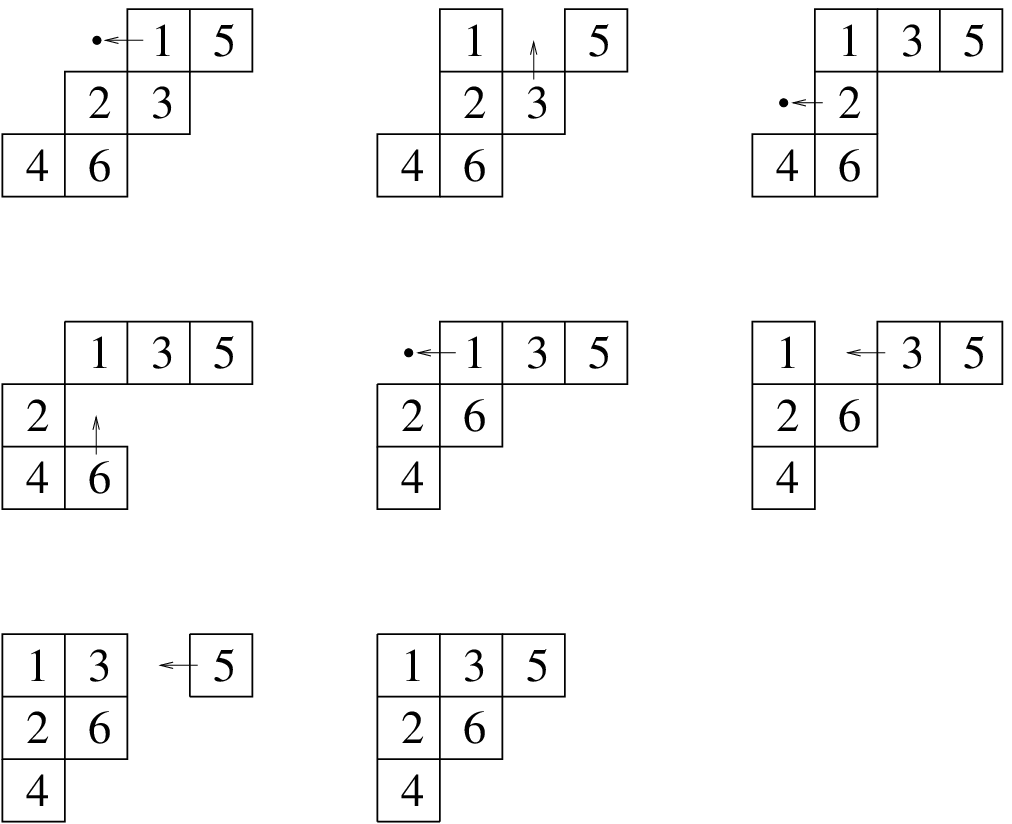,width=.50\textwidth}}
\caption{El juego de taqu\'\i n.}
\label{f:jdt}
\end{figure}

En el juego de taqu\'\i n, queremos trasformar un tableau $T$ de
forma sezgada $\nu/\mu$ en un tableau corriente de $|\nu| - |\mu|$ cajas.
Puede que la mejor manera de explicar este juego sea con un ejemplo. 
La Figura~\ref{f:jdt}
muestra el juego de taqu\'\i n en el tableau
\begin{equation*}
\NantelYoung{&&\Carre{1}&\Carre{5}\cr
           &\Carre{2}&\Carre{3}\cr
           \Carre{4}&\Carre{6}\cr
}
\end{equation*}
de la forma sezgada $(4,3,2)/(2,1)$.

En el juego de taqu\'{\i}n, movemos las cajas de $T$ hacia arriba y a la izquierda seg\'un las reglas siguientes.
\begin{enumerate}[]
\item Mientras la forma del tableau est\'e sezgada,
\begin{enumerate}[(1)]
\item Elejimos cualquier cuadrado en la frontera sudeste de la forma vac\'{\i}a.
(En la Figura~\ref{f:jdt} lo hemos indicado con un punto.)
\item Movemos a esta posici\'on la caja que se encuentra debajo o a su derecha; de estas dos, escogemos la que tenga el menor n\'umero. 
\item Mientras el tableau no tenga forma ni corriente ni sezgada, 
movemos a la posici\'on m\'as recientemente 
desocupada la caja que se encuentra debajo o a su derecha; de estas dos, escogemos la que tenga el menor n\'umero.
\end{enumerate}
\end{enumerate}

Llamamos el tableau que resulta $jdt(T)$.  Sorprendentemente, este tableau
no depende de los cuadrados elegidos sucesivamente en la regla (1).
En la expansi\'on de Schur de $s_{\nu/\mu}$ (\ref{eq:jdt}),
el coeficiente de Littlewood-Richardson 
$c_{\lambda\mu}^\nu$ es igual al n\'umero de tableaux est\'andar de Young
$T$ cuya forma es $\nu/\mu$ tales que 
$jdt(T)$ tiene la forma $\lambda$ y contiene los n\'umeros
$1,\dotsc,|\lambda|$ en ese orden, al leerlos en el \emph{orden de lectura} (de arriba a abajo y de 
izquierda a derecha).

Por ejemplo en el c\'alculo de la expansi\'on de Schur de $s_{321/21}$, 
jugamos el juego de taqu\'\i n en cada uno de los seis tableaux
est\'andar de Young de la forma $(3,2,1)/(2,1)$,
\begin{gather*}
jdt
\begin{pmatrix}
\;\;\;\;\;1\; \\
\;\;2\;\;\; \\
3\;\;\;\;\;
\end{pmatrix}
= 
\begin{matrix}
1 \\
2 \\
3 \\
\end{matrix}, \quad
jdt
\begin{pmatrix}
\;\;\;\;\; 1 \\
\;\;3\;\; \\
2\;\;\;\;
\end{pmatrix}
=
\begin{matrix}
13 \\
2\;\;\; 
\end{matrix}, \quad
jdt
\begin{pmatrix}
\;\;\;\;2\; \\
\;\;1\;\;\; \\
3\;\;\;\;\;
\end{pmatrix}
=
\begin{matrix}
12 \\
3\;\;\; 
\end{matrix}, \\
jdt
\begin{pmatrix}
\;\;\;\;2 \\
\;3\;\; \\
1\;\;\;\;\;
\end{pmatrix}
=
\begin{matrix}
12 \\
3\;\;\; 
\end{matrix}, \quad
jdt
\begin{pmatrix}
\;\;\;3\; \\
\;1\;\; \\
2\;\;\;\;
\end{pmatrix}
=
\begin{matrix}
13 \\
2\;\;\; 
\end{matrix}, \quad
jdt
\begin{pmatrix}
\;\;\;3\; \\
\;2\;\;\; \\
1\;\;\;\;\;
\end{pmatrix}
=
\begin{matrix}
123
\end{matrix}.
\end{gather*}
Ignorando el segundo y el quinto tableaux,
que no contienen los n\'umeros $1,2,3$ en el orden de lectura,
tenemos que
\begin{equation*}
s_{321/21} = s_{111} + 2s_{21} + s_3.
\end{equation*}

Como ya hemos notado (\ref{eq:tos}), no todas las funciones monomio-positivas
son Schur-positivas.  Las funciones Schur-positivas son interesantes 
porque corresponden a representaciones de $S_n$.  
(Vea por ejemplo la segunda parte de esta serie \cite{ALRS2}.)
Hay una gran cantidad de resultados y conjeturas que dicen que ciertas funciones
sim\'etricas son Schur-positivas.
(Vea, por ejemplo~\cite{StanAG}, \cite{StanPos}.)

Los resultados de la Secci\'on~\ref{s:ptnn}
se pueden aplicar a las matrices de Jacobi-Trudi
para definir funciones que son monomio-positivas.
Por ejemplo, el polinomio totalmente no negativo
$\Delta_{13,13}\Delta_{2,2} - \Delta_{12,12}\Delta_{3,3}$
nos da una clase infinita de diferencias de productos de funciones de Schur.
Curiosamente, estos polinomios parecen ser siempre Schur-positivos.
Por ejemplo, tenemos que
\begin{equation}\label{eq:schurpos}
\begin{split}
s_{41/1}s_2 - s_{32}s_1 &= s_6 + 2 s_{51} + s_{42} + s_{411}, \\
s_{52/1}s_3 - s_{43}s_2 &= s_{81} + 2s_{72} + s_{711} +  s_{63} + 2s_{621} + s_{54} + s_{531} + s_{522}, \\
s_{63/1}s_4 - s_{54}s_3 &= s_{10,2} + 2s_{93} + s_{921} + s_{84} + s_{831} 
+ s_{822} + s_{75} + s_{741}
+ 2s_{732} \\
& \quad+ s_{66} + s_{651} + s_{642} + s_{633},\\
&\vdots
\end{split}
\end{equation}

Varios casos especiales se han demostrado ya.
Por ejemplo, podemos mostrar que todos las funciones sim\'etricas de
(\ref{eq:schurpos}) son Schur-positivas.

\begin{prop}\label{p:spos}
Sea $i$ un entero no negativo.  Entonces la funci\'on sim\'etrica
\begin{equation}\label{eq:sdiff}
s_{i+3,i/1}s_{i+1} - s_{i+2,i+1}s_i
\end{equation}
es Schur-positiva.
\end{prop}
\begin{dem}
Usando el juego de taqu\'\i n y  dos interpretaciones de los coeficientes de
Littlewood-Richardson (\ref{eq:lr}), (\ref{eq:jdt}), se puede mostrar que 
tenemos que
\begin{equation*}
s_{i+3,i/1} = s_{i+2,i} + s_{i+3,i-1}.
\end{equation*}
Entonces nuestra funci\'on simétrica (\ref{eq:sdiff}) se puede escribir como
\begin{equation}\label{eq:sdiff2}
s_{i+2,i}s_{i+1} + s_{i+3,i-1}s_{i+1} - s_{i+2,i+1}s_i.
\end{equation}
Ahora tenemos que mostrar que cualquier funci\'on de Schur que aparece
en la expansión de Schur de $s_{i+2,i+1}s_i$
tambi\'en aparece en la expansi\'on de Schur de 
$s_{i+2,i}s_{i+1}$ o de $s_{i+3,i-1}s_{i+1}$.
Supongamos que aplicamos nuestro primer m\'etodo de Littlewood-Richardson al 
producto $s_{i+2,i+1}s_i$ y
consideramos una manera de rellenar el tableau de la forma $i$ y
a\~nadir su contenido a la partici\'on
$(i+2,i+1)$.  Sea $\nu$ la partici\'on que resulta. 
Aplicando ahora el m\'etodo de Littlewood-Richardson al
producto $s_{i+2,i}s_{i+1}$, 
rellenamos un tableau de la forma $i+1$ 
con los mismos n\'umeros que hemos usado para rellenar el tableau de
la forma $i$, m\'as un $2$ adicional.  A\~nadiendo el contenido de 
este tableau a la partici\'on $(i+2,i)$, obtenemos otra vez la partici\'on
$\nu$.  

Ya que $s_\nu$ aparece en la expansi\'on de Schur de 
$s_{i+2,i}s_{i+1}$ con multiplicidad mayor que o igual a su
multiplicidad en la expansi\'on de Schur de $s_{i+2,i+1}s_i$, sabemos que
la funci\'on
\begin{equation*}
s_{i+2,i}s_{i+1} - s_{i+2,i+1}s_i.
\end{equation*}
es Schur-positiva.
Entonces la funci\'on (\ref{eq:sdiff2}) tambi\'en 
es Schur-positiva.
\end{dem}  

La generalizaci\'on de la
Proposici\'on~\ref{p:spos} tambi\'en es cierta~\cite{RSkanKLImm}.

\begin{thm}\label{t:SNNpoly}
Sea $f = \Delta_{J,J'}\Delta_{L,L'} - \Delta_{I,I'}\Delta_{K,K'}$ 
un polinomio totalmente no negativo y sea $A$ una matriz de Jacobi-Trudi. 
Entonces $f(A)$ es una funci\'on Schur-positiva.
\end{thm}
Es natural preguntarse  si la no negatividad total de un polinomio $f$
implica que $f(A)$ es una funci\'on Schur-positiva por cada matriz $A$
de Jacobi-Trudi.  El cuarto autor y Andrei Zelevinsky 
han mostrado que esto no es cierto \cite{SkanNNDCB}.
Así tenemos la pregunta siguiente.

\begin{question}
?`Existe una caracterizaci\'on sencilla de todos los polinomios $f$
con la propiedad de que para cada matriz $A$ de Jacobi-Trudi,
$f(A)$ es una funci\'on Schur-positiva?
\end{question}

El Teorema~\ref{t:SNNpoly} se mostr\'o en \cite{RSkanKLImm}, utilizando
resultados de Kazhdan-Lusztig~\cite{KLRepCH} y de
Haiman~\cite{HaimanHecke} sobre las representaciones del grupo sim\'etrico
y sus caracteres.  Ya veremos m\'as de este tema en la segunda parte 
de nuestra serie de art\'\i culos. \cite{ALRS2}

\section*{\textsf{Agradecimientos}\label{s:agrad}}
Los autores quiren agradecer a
Carlos Montenegro por su ayuda en la organizacion del Primer Encuentro Colombiano de Combinatoria, y a Jason Asher, Brian Drake, Sergey Fomin, Richard Stanley, y John Stembridge por varias discusiones matem\'aticas que nos fueron de gran utilidad.

\bibliographystyle{amsplain}

\bibliography{my}

\end{document}